\documentclass[a4paper]{amsart}
\usepackage{amssymb}
\usepackage{mathtools}
\usepackage{microtype}
\usepackage{enumitem}
\usepackage{tikz}
\usetikzlibrary{calc,cd}
\usepackage{hyperref}
\hypersetup{
    colorlinks=true,
    linkcolor=purple,
    urlcolor=orange,
    citecolor=teal,
    pdfauthor = {Kanta Koeda, Ryuma Orita, and Kanon Yashiro},
    pdftitle = {Floer-type bipersistence modules and rectangle barcodes},
    pdfpagemode = UseNone,
    bookmarksnumbered=true,
}

\title{Floer-type bipersistence modules and rectangle barcodes}

\author{Kanta Koeda}
\address[Kanta Koeda]{}
\email{\href{mailto:kanta123019990513@gmail.com}{kanta123019990513@gmail.com}}

\author{Ryuma Orita}
\address[Ryuma Orita]{Department of Mathematics, Faculty of Science, Niigata University, Niigata 950-2181, Japan}
\email{\href{mailto:orita@math.sc.niigata-u.ac.jp}{orita@math.sc.niigata-u.ac.jp}}
\urladdr{\url{https://ryuma-orita.netlify.app/}}

\author{Kanon Yashiro}
\address[Kanon Yashiro]{Graduate School of Science and Technology, Niigata University, Niigata 950-2181, Japan}
\email{\href{mailto:yashiro@m.sc.niigata-u.ac.jp}{yashiro@m.sc.niigata-u.ac.jp}}
\urladdr{\url{https://sites.google.com/view/kanon-yashiro/}}

\thanks{R.~Orita was supported by JSPS KAKENHI Grant Number 21K13787.}

\newtheorem{theorem}{Theorem}[section]
\newtheorem{lemma}[theorem]{Lemma}
\newtheorem{proposition}[theorem]{Proposition}
\newtheorem{corollary}[theorem]{Corollary}

\theoremstyle{definition}
\newtheorem{definition}[theorem]{Definition}
\newtheorem{example}[theorem]{Example}

\theoremstyle{remark}
\newtheorem{remark}[theorem]{Remark}

\newcommand{\ZZ}{\mathbb{Z}}
\newcommand{\RR}{\mathbb{R}}
\newcommand{\Ztwo}{\mathbb{Z}/2\mathbb{Z}}
\DeclareMathOperator{\Image}{Im}
\setlist[enumerate]{label=(\roman*)}

\providecommand*{\given}{}
\DeclarePairedDelimiterX{\Set}[1]{\lbrace}{\rbrace}{%
  \nonscript\,\renewcommand*{\given}{\nonscript\;\delimsize\vert\nonscript\;}#1\nonscript\,}

\DeclareMathOperator{\osc}{osc}
\newcommand{\Ham}{\mathrm{Ham}}
\newcommand{\LL}{\mathcal{L}}
\newcommand{\Cr}{\mathrm{Crit}}
\DeclareMathOperator{\Coim}{Coim}
\newcommand{\FF}{\mathbb{F}}
\newcommand{\HH}{\mathbb{H}}
\newcommand{\CM}{\mathrm{CM}}
\newcommand{\HM}{\mathrm{HM}}
\newcommand{\CF}{\mathrm{CF}}
\newcommand{\HF}{\mathrm{HF}}
\newcommand{\HHMM}{\mathbb{HM}}
\newcommand{\HHFF}{\mathbb{HF}}

\begin{document}

\begin{abstract}
    We study the bipersistence modules obtained from the action-window homology of Floer-type complexes over $\Lambda^{\FF,\{0\}}=\FF$. We prove that these modules are rectangle-decomposable and establish an explicit dictionary between their graded rectangle barcodes and the classical one-parameter barcodes. This dictionary is an isometry for the bottleneck distance. Consequently, the interleaving and bottleneck distances coincide on this class of bipersistence modules: the optimal constant in algebraic stability improves from $3$ for general rectangle-decomposable $\RR^2$-indexed modules to $1$, even though the rectangles that occur are not blocks. The resulting rectangle barcodes depend $1$-Lipschitz-continuously on the Morse function in the uniform norm and on the Hamiltonian diffeomorphism in Hofer's metric. Their corner multiplicities record critical points and one-periodic Hamiltonian orbits degree by degree and filtration level by filtration level, while their geometry displays spectral invariants, boundary depth, and spectral spreads. Thus the rectangle barcode gives a geometric presentation, rather than a refinement, of the classical barcode.
\end{abstract}

\maketitle

\section{Introduction}\label{sec:introduction}

\subsection{Background}

Persistent homology is a powerful tool in Topological Data Analysis for exploring the \textit{shape} of data (see \cite{Ca09,Gh08} for surveys).
Given a partially ordered set $P$ viewed as a category, a \textit{$P$-indexed persistence module} is a functor from $P$ to the category of vector spaces over a field $\FF$ (Definition \ref{def:persistence_module}).
When $P$ is totally ordered, the Normal Form Theorem (Theorem \ref{thm:normal_form}) decomposes every pointwise finite-dimensional (p.f.d.)\ persistence module into interval modules, which leads to the notion of the \textit{barcode}. Already for $P=\RR\times\RR$, the setting of \textit{bipersistence} (or two-parameter persistence) modules, no such decomposition exists in general (see, e.g., \cite[Example 4.6]{BL23}).
A bipersistence module is called \textit{rectangle-decomposable} if it decomposes into interval modules whose supporting intervals are rectangles $I_1\times I_2$ (Definition \ref{def:rectangle-decomposable}).

\subsection{Results}\label{sec:intro_results}

Usher and Zhang \cite{UZ16} introduced Floer-type complexes $(C_*,\partial,\ell)$ over a Novikov field $\Lambda^{\FF,\Gamma}$ (see Section \ref{sec:Floer-type_complex}) and studied the associated one-parameter persistence modules $\HH_k^{\leq\bullet}(C_*)$ defined by the sublevel filtration $C_*^{\leq t}=\{\,x\in C_*\mid \ell(x)\leq t\,\}$.
In this paper, we investigate the two-parameter persistence modules $\HH_k^{(\bullet,\bullet]}(C_*)$ defined by the action windows $C_*^{(a,b]}=C_*^{\leq b}/C_*^{\leq a}$ (see Section \ref{sec:Floer-type_complex}).
Our first main theorem is the following.

\begin{theorem}[Rectangle Barcode Theorem]\label{thm:main}
    Let $(C_*,\partial,\ell)$ be a Floer-type complex over $\Lambda^{\FF,\{0\}}=\FF$ and let $k\in\ZZ$.
    Then there exists a multiset $\mathcal{B}_k^{(\bullet,\bullet]}(C_*)$ of rectangles in $\RR^2$ such that
    \[
        \HH_k^{(\bullet,\bullet]}(C_*) \cong \bigoplus_{R\in\mathcal{B}_k^{(\bullet,\bullet]}(C_*)} \FF_R.
    \]
    Moreover, the multiset $\mathcal{B}_k^{(\bullet,\bullet]}(C_*)$ is uniquely determined by $(C_*,\partial,\ell)$, and every rectangle in it is of type $(\mathrm{S})$, $(\mathrm{B})$, or $(\mathrm{N})$ described below, with bottom-right corner on the diagonal $\Delta=\Set{(a,b)\in\RR^2 \given a=b}$.
\end{theorem}

Here $\FF_I$ denotes the interval module of an interval $I$ (see Example \ref{example:interval_module}).
Given $c\in\RR$ and $\ell_1,\ell_2\in (0,\infty]$,
we set $R(c;\ell_1,\ell_2)$ to be the rectangle
\[
    R(c;\ell_1,\ell_2) = [c-\ell_1,c)\times [c,c+\ell_2)
\]
whose bottom-right corner lies on the diagonal $\Delta$ of $\RR^2$. We refer to $\ell_1$ and $\ell_2$ as its horizontal and vertical side lengths, respectively.
(For $\ell_1=\infty$ we read $[c-\ell_1,c)=(-\infty,c)$, and similarly for $\ell_2=\infty$.)
A rectangle $R(c;\ell_1,\ell_2)$ is of type $(\mathrm{S})$ if $\ell_1=\ell_2=\infty$, of type $(\mathrm{B})$ if $\ell_1=\infty>\ell_2$, of type $(\mathrm{N})$ if $\ell_2=\infty>\ell_1$, and of type $(\mathrm{R})$ if both $\ell_1$ and $\ell_2$ are finite (see Figure \ref{fig:rectangle_type} in Section \ref{sec:main_results}).

Our second main theorem determines the rectangle barcode completely in terms of the barcodes $\mathcal{B}_k^{\leq\bullet}(C_*)$ of the classical one-parameter persistence modules $\HH_k^{\leq\bullet}(C_*)$.
Recall from Lemma \ref{lem:elementary} and Theorem \ref{thm:chain_normal_form}, and also from \cite[Theorem 6.2]{UZ16}, that over $\Lambda^{\FF,\{0\}}$, every bar in $\mathcal{B}_k^{\leq\bullet}(C_*)$ has the form $[c,\infty)$ or $[a,b)$ with $a<b$.

\begin{theorem}[Dictionary Theorem]\label{thm:dictionary}
    In the setting of Theorem \textup{\ref{thm:main}}, for every degree $k\in\ZZ$ we have the equality of multisets
    \begin{align*}
        \mathcal{B}_k^{(\bullet,\bullet]}(C_*)
        ={} & \Set*{R(c;\infty,\infty) \given [c,\infty)\in\mathcal{B}_k^{\leq\bullet}(C_*)} \\
        & \sqcup \Set*{R(a;\infty,b-a) \given [a,b)\in\mathcal{B}_k^{\leq\bullet}(C_*)} \\
        & \sqcup \Set*{R(b;b-a,\infty) \given [a,b)\in\mathcal{B}_{k-1}^{\leq\bullet}(C_*)},
    \end{align*}
    where each family is counted with the multiplicity of the corresponding bars.
    In particular, the two-parameter barcodes $\bigl(\mathcal{B}_k^{(\bullet,\bullet]}(C_*)\bigr)_{k\in\ZZ}$ and the one-parameter barcodes $\bigl(\mathcal{B}_k^{\leq\bullet}(C_*)\bigr)_{k\in\ZZ}$ determine each other.
\end{theorem}

Note the degree bookkeeping in Theorem \ref{thm:dictionary}.
A finite bar $[a,b)\in\mathcal{B}_k^{\leq\bullet}(C_*)$ contributes twice: it gives the type (B) rectangle $R(a;\infty,b-a)$ in degree $k$, recording the filtered cycle born at level $a$, and the type (N) rectangle $R(b;b-a,\infty)$ in degree $k+1$, recording the non-cycle that kills this cycle at level $b$.
An infinite bar $[c,\infty)\in\mathcal{B}_k^{\leq\bullet}(C_*)$ contributes the single type (S) rectangle $R(c;\infty,\infty)$ in degree $k$.
In particular, no rectangle of type (R) occurs.

Theorem \ref{thm:dictionary} shows that the rectangle barcode is \textit{not} a stronger invariant than the classical barcode: over $\Lambda^{\FF,\{0\}}$ the two contain exactly the same information.
The value of the two-parameter picture is rather that of a \textit{normal form}: the multiplicity of summands with a fixed corner records the multiplicity of generators with the corresponding filtration value and degree (Theorems \ref{thm:morse_normal_form}, \ref{thm:floer_contractible_normal_form}, and \ref{thm:floer_non-contractible_normal_form}), and familiar invariants become visible geometric features of the barcode (see Section \ref{sec:applications}, in particular Corollary \ref{cor:noncycle_depth}). When several generators have the same filtration value and degree, indexing the corresponding rectangle summands by the individual generators is in general non-canonical.

Our third main theorem states that the correspondence of Theorem \ref{thm:dictionary} is an isometry.
For graded barcodes we write $d_{\mathrm{bot}}\bigl(\mathcal{B}_*,\mathcal{B}'_*\bigr)=\sup_{k\in\ZZ} d_{\mathrm{bot}}\bigl(\mathcal{B}_k,\mathcal{B}'_k\bigr)$ (see Section \ref{sec:persistence_modules}).

\begin{theorem}[Isometry Theorem]\label{thm:isometry}
    Let $(C_*,\partial_C,\ell_C)$ and $(D_*,\partial_D,\ell_D)$ be Floer-type complexes over $\Lambda^{\FF,\{0\}}=\FF$.
    Then
    \[
        d_{\mathrm{bot}}\left(\mathcal{B}_*^{(\bullet,\bullet]}(C_*),\mathcal{B}_*^{(\bullet,\bullet]}(D_*)\right)
        = d_{\mathrm{bot}}\left(\mathcal{B}_*^{\leq\bullet}(C_*),\mathcal{B}_*^{\leq\bullet}(D_*)\right).
    \]
\end{theorem}

\begin{corollary}\label{cor:distances}
    In the setting of Theorem \textup{\ref{thm:isometry}}, all of the following quantities coincide:
    \begin{align*}
        \sup_{k\in\ZZ}\, d_{\mathrm{int}}\left(\HH_k^{(\bullet,\bullet]}(C_*),\HH_k^{(\bullet,\bullet]}(D_*)\right)
        & = d_{\mathrm{bot}}\left(\mathcal{B}_*^{(\bullet,\bullet]}(C_*),\mathcal{B}_*^{(\bullet,\bullet]}(D_*)\right) \\
        & = d_{\mathrm{bot}}\left(\mathcal{B}_*^{\leq\bullet}(C_*),\mathcal{B}_*^{\leq\bullet}(D_*)\right) \\
        & = \sup_{k\in\ZZ}\, d_{\mathrm{int}}\left(\HH_k^{\leq\bullet}(C_*),\HH_k^{\leq\bullet}(D_*)\right).
    \end{align*}
    In particular, the isometry theorem $d_{\mathrm{int}}=d_{\mathrm{bot}}$ holds for Floer-type bipersistence modules over $\Lambda^{\FF,\{0\}}$, in the graded sense displayed above.
\end{corollary}

This is in contrast with the general two-parameter situation, where the optimal bound is $d_{\mathrm{bot}}\leq 3\,d_{\mathrm{int}}$ for rectangle-decomposable modules by a theorem of Bjerkevik \cite{Bj21} (see Theorem \ref{thm:algebraic_stability} and Remark \ref{rem:constants}).
As a consequence, we obtain stability results with constant $1$; in the functional case below, this constant is sharp.
Given a function $F$ on a closed manifold $M$, let $\|F\|=\max_M |F|$ denote the uniform norm.

\begin{theorem}[Functional Stability Theorem]\label{thm:functional_stability}
    Let $M$ be a connected closed manifold.
    For any pair of Morse functions $f,g$ on $M$, we have
    \[
        d_{\mathrm{bot}}\left(\mathcal{B}_*^{(\bullet,\bullet]}(f),\mathcal{B}_*^{(\bullet,\bullet]}(g)\right) \leq \|f-g\|,
    \]
    where $\mathcal{B}_*^{(\bullet,\bullet]}(f)$ denotes the rectangle barcode of the Morse bipersistence module of $f$ (Section \ref{sec:applications_Morse}).
    The constant $1$ is optimal, as is seen by taking $g=f+c$ for a constant $c$.
\end{theorem}

\begin{theorem}[Dynamical Stability Theorem]\label{thm:dynamical_stability}
    Let $(M,\omega)$ be a connected closed symplectically aspherical symplectic manifold.
    For any pair of non-degenerate Hamiltonian diffeomorphisms $\phi,\psi\in\Ham(M)$, we have
    \[
        d_{\mathrm{bot}}\left(\mathcal{B}_*^{(\bullet,\bullet]}(\phi),\mathcal{B}_*^{(\bullet,\bullet]}(\psi)\right) \leq d_{\mathrm{Hofer}}(\phi,\psi),
    \]
    where $\mathcal{B}_*^{(\bullet,\bullet]}(\phi)$ denotes the rectangle barcode of the Hamiltonian bipersistence module of $\phi$ (Section \ref{sec:filtered_Floer_contractible}).
    Namely, the correspondence sending a Hamiltonian diffeomorphism to its rectangle barcode is $1$-Lipschitz with respect to the Hofer metric.
\end{theorem}

For non-contractible orbits on symplectically atoroidal manifolds, the corresponding Hamiltonian-level stability estimate is stated in Theorem \ref{thm:dynamical_stability_non-contractible}.
These results improve the constant $3$ obtained in an earlier version of this paper (arXiv:2312.07847v4) to $1$.

Finally, Theorem \ref{thm:dictionary} makes classical invariants directly visible in the rectangle barcode.
The image of the spectral invariant of a Hamiltonian diffeomorphism coincides with the set of real numbers $c$ for which $(c,c)$ is the corner of a rectangle of type (S), and Usher's boundary depth is the maximum of $0$ and the finite vertical side lengths of the rectangles of type (B) (see Section \ref{sec:filtered_Floer_contractible}).
Moreover, the \textit{non-cycle depth}, defined as the maximum of $0$ and the finite horizontal side lengths of the rectangles of type (N), was introduced in that earlier version as a potentially new invariant, but it turns out to coincide with the boundary depth (Corollary \ref{cor:noncycle_depth}), answering the question raised there.
For a non-trivial free homotopy class of loops, the total Floer homology vanishes (Remark \ref{rem:non-contractible_vanishing}), so that all information is carried by the action filtration; this makes the action-window modules particularly natural for non-contractible orbits. In this setting, the spectral spreads inspired by Polterovich--Shelukhin \cite{PS16} can be read off from the rectangle barcode (Proposition \ref{prop:Floer_spectral_spread}), yielding an existence criterion for non-contractible periodic orbits (Proposition \ref{prop:existence_of_non-contractible_orbits}).
The relation with G\"{u}rel's theorem \cite{Gu13} is discussed in Remark \ref{rem:Gurel}.

\subsection*{Attribution}
Theorems \ref{thm:dictionary} and \ref{thm:isometry}, together with their proofs and the reduction of Theorem \ref{thm:main} to the chain-level normal form in Section \ref{sec:main_results}, are due to the anonymous referee of an earlier version of this paper and are included here with the referee's generous permission. The responsibility for any inaccuracy in the present write-up lies solely with the authors.
Theorem \ref{thm:dictionary} was partly anticipated in 2024 by Q.~Feng and V.~Lebovici in private communications: Feng derived an equivalent description of the rectangle barcode from the singular value decomposition of \cite{UZ16}, while Lebovici observed that only the types (S), (B), and (N) occur, that the graded one- and two-parameter barcodes carry the same information, and that the Morse-theoretic case of Theorem \ref{thm:main} can be deduced from homology pyramids (Remark \ref{rem:pyramid}).
That earlier version proved Theorem \ref{thm:main} by a different multiparameter argument, using weak exactness on all of $\RR^2$ and middle exactness on the restriction to $\Set{(a,b)\in\RR^2 \given a\leq b}$ \cite{BLO22,BLO23}.
Middle exactness was asserted only on that restriction.
Indeed, the modules $\HH_k^{(\bullet,\bullet]}(C_*)$ themselves are in general not middle-exact: by the theorem of Cochoy and Oudot \cite{CO20} they would then be block-decomposable, whereas the occurring rectangles are not blocks.
The present chain-level proof follows the referee's argument and additionally yields the explicit dictionary of Theorem \ref{thm:dictionary}.
Worked examples of rectangle barcodes of Morse functions can be found in that version (Section 5.3).

\subsection*{Related work and outlook}
Usher and Zhang \cite{UZ16} proved that the one-parameter barcode determines the filtered chain homotopy type of a Floer-type complex and, for general $\Gamma$, introduced concise and verbose barcodes via a non-Archimedean singular value decomposition. Over $\Lambda^{\FF,\{0\}}$, Feng used this decomposition to obtain a description equivalent to Theorem \ref{thm:dictionary} (see the Attribution subsection). In classical topology, interlevel persistence of a suitably tame function carries the same information as extended persistence, through an explicit correspondence with a homological degree shift \cite{CdSM09,CEH09}.
Remark \ref{rem:pyramid} explains the corresponding homology-pyramid picture here.
Usher \cite{Ush23} developed filtered matched pairs and Poincar\'{e}--Novikov structures producing interlevel-type barcodes in Hamiltonian Floer theory, and \cite{Ush25} recast this structure in terms of filtered cospans, discussed in Remark \ref{rem:cospan}. Chach\'{o}lski, Giunti, and Landi \cite[Theorem 4.2]{CGL21} obtained a decomposition theorem for parametrized chain complexes under hypotheses that do not cover all Floer-type complexes. Polterovich and Shelukhin \cite[Definition 2.1]{PS16} considered the closely related notion of restricted two-parametric persistence modules.
Theorem \ref{thm:dictionary} shows that the bipersistence modules studied here carry exactly the same information as the classical barcode. This points to two natural ways in which genuinely two-parameter information can arise: enriching the structure by a descending filtration, as in \cite{Ush23,Ush25}, or relaxing the chain-complex condition itself. In such enlarged settings the recovery mechanism of Theorem \ref{thm:dictionary} need not persist, and invariants such as the supports of Betti numbers may contain genuinely two-parameter information.
The second route is pursued in forthcoming work of the second and third authors.

\subsection*{Organization of the paper}
Section \ref{sec:multiparameter} collects the background on persistence modules and Floer-type complexes.
Section \ref{sec:main_results} proves Theorems \ref{thm:main} and \ref{thm:dictionary}.
Section \ref{sec:stability} proves Theorem \ref{thm:isometry} and Corollary \ref{cor:distances}.
Section \ref{sec:applications} presents the applications, including the proofs of Theorems \ref{thm:functional_stability} and \ref{thm:dynamical_stability}.

\section{Preliminaries}\label{sec:multiparameter}

\subsection{Persistence modules, barcodes and distances}\label{sec:persistence_modules}

We briefly recall the necessary notions, mainly following \cite{BL23}.
Let us fix a field $\FF$ and let $(P,\preceq)$ be a partially ordered set.

\begin{definition}\label{def:interval}
    A non-empty subset $I$ of $P$ is called an \textit{interval} if $I$ satisfies the following conditions.
    \begin{enumerate}
        \item For any $x,y\in I$ and $z\in P$, the condition $x\preceq z\preceq y$ implies $z\in I$.
        \item For any $x,y\in I$ there exist $z_1,\dots,z_k\in I$ such that $z_i$ and $z_{i+1}$ are comparable for all $0\leq i\leq k$, where $z_0=x$ and $z_{k+1}=y$.
    \end{enumerate}
\end{definition}

A multiset of intervals in $P$ is called a \textit{barcode} over $P$.
For the totally ordered set $\RR$, intervals are of the forms $(x,y)$, $(x,y]$, $[x,y)$ and $[x,y]$ with $x,y\in\RR\cup\{\pm\infty\}$, with the evident conventions at infinite endpoints.
Given totally ordered sets $T_1,\dots,T_n$, we equip the product $T_1\times\dots\times T_n$ with the componentwise partial order.
A \textit{rectangle} in $T_1\times\dots\times T_n$ is an interval of the form $I_1\times\dots\times I_n$, where each $I_i$ is an interval in $T_i$.
A \textit{block} in $T_1\times T_2$ is an interval of one of the forms $D_1\times D_2$, $U_1\times U_2$, $I_1\times T_2$, or $T_1\times I_2$, where $D_i$ and $U_i$ are downsets and upsets, respectively, and the $I_i$ are intervals.
Every block is a rectangle, whereas the rectangles of types (S), (B), (N) and (R) of Section \ref{sec:intro_results} are not blocks.

\begin{definition}\label{def:persistence_module}
    A \textit{$P$-indexed persistence module} is a functor $\mathfrak{M}\colon P\to\mathbf{Vect}_\FF$, where $P$ is viewed as a category.
    We write $\mathfrak{M}_x$ for the value at $x\in P$ and $\mathfrak{M}_{x,y}\colon \mathfrak{M}_x\to \mathfrak{M}_y$ for the structure map associated to $x\preceq y$.
    A persistence module $\mathfrak{M}$ is \textit{pointwise finite-dimensional} (p.f.d.) if each $\mathfrak{M}_x$ is finite-dimensional.
\end{definition}

Direct sums of persistence modules are defined pointwise, and $\mathfrak{M}$ is \textit{indecomposable} if it admits no decomposition into two non-zero summands.

\begin{example}\label{example:interval_module}
    Let $I$ be an interval in $P$.
    The \textit{interval module} $\FF_I\colon P\to\mathbf{Vect}_\FF$ is defined by $(\FF_I)_x=\FF$ for $x\in I$ and $(\FF_I)_x=0$ otherwise, with the structure map $(\FF_I)_{x,y}=\mathrm{id}_{\FF}$ for $x\preceq y$ in $I$ and zero otherwise.
    Interval modules are indecomposable \cite[Proposition 2.2]{BL18}.
\end{example}

\begin{theorem}[{Normal Form Theorem; \cite[Theorem 1.2]{BCB20}, \cite{CB15}}]\label{thm:normal_form}
    Let $T$ be a totally ordered set and $\mathfrak{M}$ a p.f.d.\ $T$-indexed persistence module.
    Then there exists a unique multiset $\mathcal{B}(\mathfrak{M})$ of intervals in $T$ such that
    $\mathfrak{M}\cong\bigoplus_{I\in\mathcal{B}(\mathfrak{M})}\FF_I$.
\end{theorem}

\begin{definition}\label{def:rectangle-decomposable}
    A $P$-indexed persistence module $\mathfrak{M}$ is \textit{interval-decomposable} if there exists a multiset $\mathcal{B}(\mathfrak{M})$ of intervals in $P$ such that $\mathfrak{M}\cong\bigoplus_{I\in\mathcal{B}(\mathfrak{M})}\FF_I$.
    An interval-decomposable $(T_1\times T_2)$-indexed persistence module is \textit{block-decomposable} if all intervals in $\mathcal{B}(\mathfrak{M})$ are blocks, and \textit{rectangle-decomposable} if they are all rectangles.
\end{definition}

Since interval modules are indecomposable, the multiset $\mathcal{B}(\mathfrak{M})$ of an interval-decomposable module is unique by the Azumaya--Krull--Remak--Schmidt theorem \cite[Theorem 1]{Az50}, and it is called the \textit{barcode} of $\mathfrak{M}$.

We next recall the standard metrics, following \cite{CdSGO16} and \cite[Section \textup{VIII}]{BL23}.
Let $\mathfrak{M}$ be an $\RR^n$-indexed persistence module.
Given $u\in [0,\infty)^n$, the \textit{$u$-shift} $\mathfrak{M}[u]$ is defined by $\mathfrak{M}[u]_x=\mathfrak{M}_{x+u}$ and $\mathfrak{M}[u]_{x,y}=\mathfrak{M}_{x+u,y+u}$, and $\Phi_{\mathfrak{M}}^u\colon \mathfrak{M}\to\mathfrak{M}[u]$ denotes the morphism with components $\mathfrak{M}_{x,x+u}$.
For $\delta\geq 0$ we write $\vec{\delta}=(\delta,\dots,\delta)$.

\begin{definition}\label{def:interleaving_distance}
    Two $\RR^n$-indexed persistence modules $\mathfrak{M}$ and $\mathfrak{N}$ are \textit{$\delta$-interleaved} if there exist morphisms $F\colon \mathfrak{M}\to\mathfrak{N}[\vec{\delta}]$ and $G\colon \mathfrak{N}\to\mathfrak{M}[\vec{\delta}]$ such that
    $G[\vec{\delta}]\circ F = \Phi_{\mathfrak{M}}^{2\vec{\delta}}$ and $F[\vec{\delta}]\circ G = \Phi_{\mathfrak{N}}^{2\vec{\delta}}$.
    The \textit{interleaving distance} is
    $d_{\mathrm{int}}(\mathfrak{M},\mathfrak{N})=\inf\Set{\delta\geq 0 \given \mathfrak{M},\mathfrak{N}\ \text{are $\delta$-interleaved}}$.
\end{definition}

A \textit{matching} $\mu\colon X\nrightarrow Y$ between multisets is a bijection $\mu\colon X'\to Y'$ between subsets $X'\subset X$ and $Y'\subset Y$.
We write $\Coim\mu=X'$ and $\Image\mu=Y'$.
Two intervals $I,J$ in $\RR^n$ are \textit{$\delta$-interleaved} if the interval modules $\FF_I,\FF_J$ are $\delta$-interleaved.
An interval $I$ is \textit{$\delta$-trivial} if $x+\vec{\delta}\notin I$ for every $x\in I$.
Given a barcode $\mathcal{B}$, let $\mathcal{B}_\delta\subset\mathcal{B}$ denote the multiset of intervals which are \textit{not} $\delta$-trivial.

\begin{definition}\label{def:bottleneck_distance}
    A \textit{$\delta$-matching} between barcodes $\mathcal{B}$ and $\mathcal{C}$ is a matching $\mu\colon\mathcal{B}\nrightarrow\mathcal{C}$ such that $\mathcal{B}_{2\delta}\subset\Coim\mu$, $\mathcal{C}_{2\delta}\subset\Image\mu$, and $I$ and $\mu(I)$ are $\delta$-interleaved for every $I\in\Coim\mu$.
    The \textit{bottleneck distance} is
    \[
        d_{\mathrm{bot}}(\mathcal{B},\mathcal{C})=\inf\Set*{\delta\geq 0 \given \text{there exists a $\delta$-matching between $\mathcal{B}$ and $\mathcal{C}$}}.
    \]
    For graded barcodes $\mathcal{B}_*=(\mathcal{B}_k)_{k\in\ZZ}$ we set $d_{\mathrm{bot}}(\mathcal{B}_*,\mathcal{C}_*)=\sup_{k\in\ZZ}d_{\mathrm{bot}}(\mathcal{B}_k,\mathcal{C}_k)$.
\end{definition}

\begin{remark}\label{rem:trivial_interleaving}
    For an interval $I$, we have $\Phi_{\FF_I}^{2\vec{\delta}}\neq 0$ if and only if $I$ is not $2\delta$-trivial.
    Consequently, any two $2\delta$-trivial intervals are $\delta$-interleaved via the zero morphisms.
\end{remark}

For p.f.d.\ one-parameter persistence modules, the interleaving and bottleneck distances coincide, which is the celebrated isometry theorem.

\begin{theorem}[\cite{CSEH07,CCSGGO09,CdSGO16,BL15,BL20,Bj21}]\label{thm:one-para_isometry}
    For any p.f.d.\ $\RR$-indexed persistence modules $\mathfrak{M}$ and $\mathfrak{N}$, we have
    \[
        d_{\mathrm{int}}(\mathfrak{M},\mathfrak{N})=d_{\mathrm{bot}}\bigl(\mathcal{B}(\mathfrak{M}),\mathcal{B}(\mathfrak{N})\bigr).
    \]
\end{theorem}

In the multiparameter setting one still has the \textit{converse algebraic stability theorem}.

\begin{proposition}[{\cite[Proposition 2.13]{BL18}}]\label{thm:converse_algebraic_stability}
    For any interval-decomposable $\RR^n$-indexed persistence modules $\mathfrak{M}$ and $\mathfrak{N}$, we have
    $d_{\mathrm{int}}(\mathfrak{M},\mathfrak{N})\leq d_{\mathrm{bot}}(\mathcal{B}(\mathfrak{M}),\mathcal{B}(\mathfrak{N}))$.
\end{proposition}

The algebraic stability theorem $d_{\mathrm{bot}}\leq d_{\mathrm{int}}$ fails in general.
The following bounds of Bjerkevik are optimal.

\begin{theorem}[{\cite[Theorems 4.3 and 4.18]{Bj21}}]\label{thm:algebraic_stability}
    Let $\mathfrak{M}$ and $\mathfrak{N}$ be p.f.d.\ rectangle-decomposable $\RR^n$-indexed persistence modules.
    Then
    \[
        d_{\mathrm{bot}}\bigl(\mathcal{B}(\mathfrak{M}),\mathcal{B}(\mathfrak{N})\bigr)\leq (2n-1)\,d_{\mathrm{int}}(\mathfrak{M},\mathfrak{N}),
    \]
    and the bound is optimal for $n=2$.
    If $n=2$ and $\mathfrak{M}$ and $\mathfrak{N}$ are block-decomposable, then the sharper bound
    \[
        d_{\mathrm{bot}}\bigl(\mathcal{B}(\mathfrak{M}),\mathcal{B}(\mathfrak{N})\bigr)\leq d_{\mathrm{int}}(\mathfrak{M},\mathfrak{N})
    \]
    holds.
\end{theorem}

\subsection{Floer-type complexes and their persistence modules}\label{sec:Floer-type_complex}

We recall the formalism of Usher and Zhang \cite{UZ16}.
Fix a ground field $\FF$ and an additive subgroup $\Gamma\subset\RR$.
The \textit{Novikov field} is
\[
    \Lambda^{\FF,\Gamma}=\Set*{\lambda=\sum_{\gamma\in\Gamma}a_\gamma s^\gamma \given a_\gamma\in\FF\ \text{and}\ \#\Set{\gamma \given a_\gamma\neq 0,\ \gamma<C}<\infty\ \text{for all}\ C\in\RR},
\]
where $s$ is a formal symbol.
Note that $\Lambda^{\FF,\{0\}}=\FF$.
It carries the valuation $\nu(0)=+\infty$ and, for $\lambda\neq 0$, $\nu(\lambda)=\min\Set{\gamma \given a_\gamma\neq 0}$.
A \textit{non-Archimedean normed vector space} over $\Lambda^{\FF,\Gamma}$ is a pair $(C,\ell)$ of a $\Lambda^{\FF,\Gamma}$-vector space $C$ and a function $\ell\colon C\to\RR\cup\{-\infty\}$ such that $\ell(x)=-\infty$ if and only if $x=0$, $\ell(\lambda x)=\ell(x)-\nu(\lambda)$, and $\ell(x+y)\leq\max\{\ell(x),\ell(y)\}$ \cite[Definition 2.2]{UZ16}.
A finite-dimensional $(C,\ell)$ is \textit{orthogonalizable} if it admits a basis $(w_1,\dots,w_n)$ with $\ell\bigl(\sum_i\lambda_iw_i\bigr)=\max_i\bigl(\ell(w_i)-\nu(\lambda_i)\bigr)$ \cite[Definition 2.10]{UZ16}.

\begin{definition}[{\cite[Definition 4.1]{UZ16}}]\label{def:Floer-type_complex}
    A \textit{Floer-type complex} $(C_*,\partial,\ell)$ over $\Lambda^{\FF,\Gamma}$ is a chain complex $(C_*=\bigoplus_{k\in\ZZ}C_k,\partial)$ over $\Lambda^{\FF,\Gamma}$ equipped with a function $\ell\colon C_*\to\RR\cup\{-\infty\}$ such that each $(C_k,\ell|_{C_k})$ is orthogonalizable and $\ell(\partial x)\leq\ell(x)$ for all $x$.
\end{definition}

Given Floer-type complexes $(C_*,\partial_C,\ell_C)$ and $(D_*,\partial_D,\ell_D)$ over $\Lambda^{\FF,\Gamma}$, their \textit{direct sum} is the Floer-type complex $C_*\oplus D_*$ with differential $\partial_C\oplus\partial_D$ and filtration $\ell(x,y)=\max\{\ell_C(x),\ell_D(y)\}$.
A chain isomorphism $\Phi\colon C_*\to D_*$ satisfying $\ell_D(\Phi x)=\ell_C(x)$ for all $x\in C_*$ is called a \textit{filtered chain isomorphism}.
If such a $\Phi$ exists, the two complexes are called \textit{filtered chain isomorphic}.

For $t\in\RR$ we set $C_k^{\leq t}=\Set{x\in C_k \given \ell(x)\leq t}$.
Since $\ell(\partial x)\leq\ell(x)$, we obtain a subcomplex $(C_*^{\leq t},\partial)$, whose homology we denote by $H_k^{\leq t}=H_k(C_*^{\leq t})$.
For $t\leq t'$ the inclusion induces $\iota_{tt'}\colon H_k^{\leq t}\to H_k^{\leq t'}$.

\begin{definition}\label{def:Floer-type_persistence}
    Given a Floer-type complex $(C_*,\partial,\ell)$ over $\Lambda^{\FF,\Gamma}$ and $k\in\ZZ$,
    the $\RR$-indexed persistence module $\bigl(\{H_k^{\leq t}\}_{t\in\RR},\{\iota_{tt'}\}_{t\leq t'}\bigr)$ is called the \textit{Floer-type persistence module in degree $k$}, denoted by $\HH_k^{\leq\bullet}(C_*)$.
\end{definition}

When $\Gamma=\{0\}$, each $C_k$ is finite-dimensional over $\FF$, so $\HH_k^{\leq\bullet}(C_*)$ is p.f.d.\ and Theorem \ref{thm:normal_form} provides its barcode $\mathcal{B}_k^{\leq\bullet}(C_*)$.
In this case we write $\mathcal{B}_*^{\leq\bullet}(C_*)=\bigl(\mathcal{B}_k^{\leq\bullet}(C_*)\bigr)_{k\in\ZZ}$ for the graded barcode. For general $\Gamma$, finite-dimensionality over the Novikov field does not in general imply pointwise finite-dimensionality over $\FF$, so we do not use the ordinary p.f.d.\ barcode formalism without additional hypotheses.

For $a\leq b$ we define the quotient complex
\[
    C_*^{(a,b]}=C_*^{\leq b}/C_*^{\leq a}
\]
and call its homology $H_k^{(a,b]}=H_k\bigl(C_*^{(a,b]}\bigr)$ the \textit{$k$-th action-window Floer-type homology}.
The interval $(a,b]$ is referred to as an \textit{action window} \cite{FH94}.\footnote{Earlier versions of this paper used the term \textit{interlevel window}. Since $H_k^{(a,b]}$ is a counterpart of the relative homology of a pair of sublevel sets rather than of the absolute homology of an interlevel set, we follow the Floer-theoretic terminology here. See the discussion in \cite[Section 1]{Ush23}.}
For $a>b$ we set $C_*^{(a,b]}=0$ and $H_*^{(a,b]}=0$.
For $b\leq b'$ the inclusion $C_k^{\leq b}\hookrightarrow C_k^{\leq b'}$ descends to the quotients and induces $\iota_{bb'}^{a}\colon H_k^{(a,b]}\to H_k^{(a,b']}$.
For $a\leq a'\leq b$ the quotient map descends to $C_k^{(a,b]}\to C_k^{(a',b]}$ and induces $\pi_{aa'}^{b}\colon H_k^{(a,b]}\to H_k^{(a',b]}$.
For $a'>b$ we set $\pi_{aa'}^{b}=0$.
Together with the unique maps whose source is a zero complex, these maps define the structure maps for all comparable pairs in $\RR^2$.
These maps commute:
\begin{equation*}
    \begin{tikzcd}
    H_k^{(a,b']} \arrow[r, "\pi_{aa'}^{b'}"] & H_k^{(a',b']} \\
    H_k^{(a,b]} \arrow[u, "\iota_{bb'}^{a}"] \arrow[r, "\pi_{aa'}^{b}"'] & H_k^{(a',b]} \arrow[u, "\iota_{bb'}^{a'}"']
    \end{tikzcd}
\end{equation*}

\begin{definition}\label{def:Floer-type_bipersistence}
    Given a Floer-type complex $(C_*,\partial,\ell)$ over $\Lambda^{\FF,\Gamma}$ and $k\in\ZZ$, the $\RR^2$-indexed persistence module
    \[
        \left(\{H_k^{(a,b]}\}_{(a,b)\in\RR^2},\{\iota_{bb'}^{a}\}_{a\in\RR,\,b\leq b'},\{\pi_{aa'}^{b}\}_{b\in\RR,\,a\leq a'}\right)
    \]
    is called the \textit{Floer-type bipersistence module in degree $k$}, denoted by $\HH_k^{(\bullet,\bullet]}(C_*)$.
\end{definition}

When $\Gamma=\{0\}$, the same finite-dimensionality argument shows that $\HH_k^{(\bullet,\bullet]}(C_*)$ is p.f.d. For general $\Gamma$, pointwise finite-dimensionality over $\FF$ need not hold.

\section{The rectangle barcode via the chain-level normal form}\label{sec:main_results}

In this section we prove Theorems \ref{thm:main} and \ref{thm:dictionary}.
Throughout this section, $\Gamma=\{0\}$, so that $\Lambda^{\FF,\{0\}}=\FF$ and $\nu(\lambda)=0$ for all $\lambda\neq 0$.

Recall from Section \ref{sec:intro_results} the rectangles $R(c;\ell_1,\ell_2)=[c-\ell_1,c)\times[c,c+\ell_2)$ with bottom-right corner on the diagonal $\Delta$.
Such rectangles can be classified into the following four types (Figure \ref{fig:rectangle_type}):
\begin{itemize}
    \item[(S)] $R(c;\infty,\infty)$, unbounded in both coordinate directions.
    \item[(B)] $R(c;\infty,\ell_2)$, unbounded horizontally and with finite vertical side length $\ell_2$.
    \item[(N)] $R(c;\ell_1,\infty)$, with finite horizontal side length $\ell_1$ and unbounded vertically.
    \item[(R)] $R(c;\ell_1,\ell_2)$, with both side lengths finite.
\end{itemize}
The initials S, B, N, and R correspond to ``spectral invariant,'' ``boundary depth,'' ``non-cycle depth,'' and ``true rectangle,'' respectively.
\begin{figure}[ht]
\centering
\begin{tikzpicture}
\draw (-2.5,-2.5) -- (3,3) node [anchor=south west] {$\Delta$};
\fill[gray,opacity=0.2] (-3.5,-1) rectangle(-1,3.5);
\fill[teal,opacity=0.2] (-3.5,1.5) rectangle(1.5,2.5);
\fill[purple,opacity=0.2] (-3,-2) rectangle(-2,3.5);
\fill[orange,opacity=0.2] (-1.5,0) rectangle(0,1);
\draw[gray,dashed,very thick] (-1,3.5) -- (-1,-1);
\draw[gray,very thick] (-3.5,-1) -- (-1,-1);
\draw[teal,dashed,very thick] (-3.5,2.5) -- (1.5,2.5) -- (1.5,1.5);
\draw[teal,very thick] (-3.5,1.5) -- (1.5,1.5);
\draw[purple,dashed,very thick] (-2,3.5) -- (-2,-2);
\draw[purple,very thick] (-3,3.5) -- (-3,-2) -- (-2,-2);
\draw[orange,dashed,very thick] (-1.5,1) -- (0,1) -- (0,0);
\draw[orange,very thick] (-1.5,1) -- (-1.5,0) -- (0,0);
\node[gray,anchor=north west] at (-1,-1) {Type (S)};
\node[teal,anchor=north west] at (1.5,1.5) {Type (B)};
\node[purple,anchor=north west] at (-2,-2) {Type (N)};
\node[orange,anchor=north west] at (0,0) {Type (R)};
\draw[gray,very thick] (-5.5,-1) -- (-5.5,3.5);
\draw[gray] (-5.5,-1.5).. controls ($(-5.5,-1.5)!.2!(-5.5,-1)!10pt!90:(-5.5,-1)$) and ($(-5.5,-1.5)!.8!(-5.5,-1)!10pt!90:(-5.5,-1)$) .. (-5.5,-1);
\fill[gray] (-5.5,-1) circle (2pt) node [fill=white,below=2mm] {spectral value $c$};
\draw[teal] (-6,1.5).. controls ($(-6,1.5)!.2!(-6,2.5)!10pt!90:(-6,2.5)$) and ($(-6,1.5)!.8!(-6,2.5)!10pt!90:(-6,2.5)$) .. (-6,2.5) node [midway,fill=white,left=-1.5mm] {$\ell_2$};
\draw[teal,very thick] (-6,1.5) -- (-6,2.5);
\fill[teal] (-6,1.5) circle (2pt);
\draw[teal,fill=white] (-6,2.5) circle (2pt);
\draw[->] (-4,1) -- (-5,1) node[midway,above] {$-\infty$};
\draw[gray,dotted] (-3.65,-1) -- (-5.35,-1);
\draw[teal,dotted] (-3.65,1.5) -- (-5.35,1.5);
\draw[teal,dotted] (-3.65,2.5) -- (-5.35,2.5);
\draw[teal,dotted] (-5.65,1.5) -- (-5.85,1.5);
\draw[teal,dotted] (-5.65,2.5) -- (-5.85,2.5);
\end{tikzpicture}
\caption{Rectangles of the four types. A type-(S) corner records a spectral value. The finite vertical side length of a type-(B) rectangle and the finite horizontal side length of a type-(N) rectangle record the corresponding bar length; their maxima give the boundary depth and the non-cycle depth, respectively. Solid boundary edges belong to the (half-open) rectangles, dashed ones do not; unbounded sides are truncated by the frame.}
\label{fig:rectangle_type}
\end{figure}

\begin{lemma}\label{lem:additivity}
    Let $(C_*,\partial_C,\ell_C)$ and $(D_*,\partial_D,\ell_D)$ be Floer-type complexes over $\Lambda^{\FF,\Gamma}$.
    \begin{enumerate}
        \item For every $k\in\ZZ$ we have $\HH_k^{(\bullet,\bullet]}(C_*\oplus D_*)\cong\HH_k^{(\bullet,\bullet]}(C_*)\oplus \HH_k^{(\bullet,\bullet]}(D_*)$ and $\HH_k^{\leq\bullet}(C_*\oplus D_*)\cong\HH_k^{\leq\bullet}(C_*)\oplus \HH_k^{\leq\bullet}(D_*)$.
        \item If $C_*$ and $D_*$ are filtered chain isomorphic, then $\HH_k^{(\bullet,\bullet]}(C_*)\cong \HH_k^{(\bullet,\bullet]}(D_*)$ and $\HH_k^{\leq\bullet}(C_*)\cong \HH_k^{\leq\bullet}(D_*)$ for every $k\in\ZZ$.
    \end{enumerate}
\end{lemma}

\begin{proof}
(i) Since $\ell(x,y)=\max\{\ell_C(x),\ell_D(y)\}$, we have $(C_*\oplus D_*)^{\leq t}=C_*^{\leq t}\oplus D_*^{\leq t}$ for every $t$, hence $(C_*\oplus D_*)^{(a,b]}=C_*^{(a,b]}\oplus D_*^{(a,b]}$ as chain complexes, compatibly with all the maps $\iota$ and $\pi$.
Taking homology yields the claim.
(ii) A filtered chain isomorphism $\Phi$ restricts to isomorphisms $C_*^{\leq t}\to D_*^{\leq t}$ for every $t$ and therefore descends to isomorphisms of all quotient complexes $C_*^{(a,b]}\to D_*^{(a,b]}$, commuting with the maps $\iota$ and $\pi$.
Taking homology yields the claim.
\end{proof}

We next compute the invariants of the two kinds of \textit{elementary} Floer-type complexes.
For $c\in\RR$ and $k\in\ZZ$, let $E(c;k)$ denote the Floer-type complex with $E(c;k)_k=\FF x$, $\ell(x)=c$, all other degrees zero, and $\partial=0$.
For $a\leq b$ and $k\in\ZZ$, let $E(a,b;k)$ denote the Floer-type complex with $E(a,b;k)_{k+1}=\FF y$, $E(a,b;k)_k=\FF x$, all other degrees zero, $\partial y=x$, $\ell(x)=a$ and $\ell(y)=b$.

\begin{lemma}\label{lem:elementary}
    The following hold.
    \begin{enumerate}
        \item In degree $j=k$ we have $\mathcal{B}_j^{\leq\bullet}\bigl(E(c;k)\bigr)=\{[c,\infty)\}$ and $\HH_j^{(\bullet,\bullet]}\bigl(E(c;k)\bigr)\cong\FF_{R(c;\infty,\infty)}$, and both vanish otherwise.
        \item Let $a<b$. Then $\mathcal{B}_j^{\leq\bullet}\bigl(E(a,b;k)\bigr)=\{[a,b)\}$ for $j=k$ and $=\emptyset$ otherwise.
        Moreover,
        \[
            \HH_j^{(\bullet,\bullet]}\bigl(E(a,b;k)\bigr)\cong
            \begin{cases}
                \FF_{R(a;\infty,b-a)} & \text{if}\ j=k,\\
                \FF_{R(b;b-a,\infty)} & \text{if}\ j=k+1,\\
                0 & \text{otherwise}.
            \end{cases}
        \]
        \item If $a=b$, then $\mathcal{B}_j^{\leq\bullet}\bigl(E(a,a;k)\bigr)=\emptyset$ and $\HH_j^{(\bullet,\bullet]}\bigl(E(a,a;k)\bigr)=0$ for every $j$.
    \end{enumerate}
\end{lemma}

\begin{proof}
(i) Since $\nu\equiv 0$ on $\FF\setminus\{0\}$, we have $\ell(\lambda x)=c$ for $\lambda\neq 0$, so $E(c;k)^{\leq t}_k=\FF x$ for $t\geq c$ and $=0$ otherwise.
As $\partial=0$, the class $[x]\in H_k^{\leq t}$ is non-zero exactly for $t\geq c$, and the structure maps send $[x]$ to $[x]$, which gives the bar $[c,\infty)$.
In the window $(a,b]$, the image $\bar{x}$ of $x$ is non-zero in $E(c;k)^{(a,b]}_k$ if and only if $\ell(x)\leq b$ and $\ell(x)>a$, that is, $(a,b)\in (-\infty,c)\times[c,\infty)=R(c;\infty,\infty)$.
All windows in this region carry the one-dimensional space $\FF\bar{x}$, and the maps $\iota$ and $\pi$ send $\bar{x}$ to $\bar{x}$, hence act as the identity between any two such windows.
This gives $\HH_k^{(\bullet,\bullet]}\cong\FF_{R(c;\infty,\infty)}$.

(ii) In the sublevel filtration, $x$ is a cycle born at $t=a$.
Since $y\in E^{\leq t}$ exactly for $t\geq b$ and $\partial y=x$, the class $[x]$ dies at $t=b$, giving the bar $[a,b)$ in degree $k$.
The generator $y$ is never a cycle in $E^{\leq t}$ (as $\partial y=x\neq 0$ in $E^{\leq t}$ whenever $y\in E^{\leq t}$), so $\mathcal{B}_{k+1}^{\leq\bullet}=\emptyset$.

In degree $k$ of the window $(a',b']$, the class $\bar{x}$ is non-zero if and only if $a'<a\leq b'$, and it lies in the image of $\bar{\partial}$ if and only if $\bar{y}$ is available, that is, $b'\geq b$.
Hence $H_k^{(a',b']}\neq 0$ exactly on $(-\infty,a)\times[a,b)=R(a;\infty,b-a)$, where it is spanned by $[\bar{x}]$, and the structure maps act as the identity.
In degree $k+1$, the class $\bar{y}$ is non-zero if and only if $a'<b\leq b'$, and $\bar{y}$ is a cycle in $E^{(a',b']}$ if and only if $\bar{\partial}\bar{y}=\bar{x}=0$, that is, $\ell(x)=a\leq a'$ (the case $b'<a$ cannot occur since $b'\geq b>a$).
Hence $H_{k+1}^{(a',b']}\neq 0$ exactly on $[a,b)\times[b,\infty)=R(b;b-a,\infty)$, spanned by $[\bar{y}]$, with identity structure maps.

(iii) If $a=b$, the regions found in (ii) are empty, and there is no bar.
\end{proof}

The following normal form theorem for Floer-type complexes over $\Lambda^{\FF,\{0\}}$ goes back to Barannikov \cite{barannikov1994}.
In the present formulation see \cite[Proposition 7.4]{UZ16}, \cite[Section 6.2]{PRSZ20} and \cite[Corollary 3.17]{Ush25}.

\begin{theorem}[Chain-level Normal Form Theorem]\label{thm:chain_normal_form}
    Every Floer-type complex $(C_*,\partial,\ell)$ over $\Lambda^{\FF,\{0\}}$ is filtered chain isomorphic to a direct sum of elementary complexes of the forms $E(c;k)$ and $E(a,b;k)$ in which only finitely many summands are non-zero in any given degree.
\end{theorem}

\begin{proof}[Proof of Theorems \textup{\ref{thm:main}} and \textup{\ref{thm:dictionary}}]
By Theorem \ref{thm:chain_normal_form}, $(C_*,\partial,\ell)$ is filtered chain isomorphic to a direct sum of elementary complexes.
By Lemma \ref{lem:additivity}, both $\HH_k^{(\bullet,\bullet]}(C_*)$ and $\HH_k^{\leq\bullet}(C_*)$, together with their barcodes, can be computed summand-wise, and by Lemma \ref{lem:elementary} each summand contributes as follows.
A summand $E(c;k)$ contributes the bar $[c,\infty)$ to $\mathcal{B}_k^{\leq\bullet}$ and the rectangle $R(c;\infty,\infty)$ to $\mathcal{B}_k^{(\bullet,\bullet]}$.
A summand $E(a,b;k)$ with $a<b$ contributes the bar $[a,b)$ to $\mathcal{B}_k^{\leq\bullet}$, the rectangle $R(a;\infty,b-a)$ to $\mathcal{B}_k^{(\bullet,\bullet]}$ and the rectangle $R(b;b-a,\infty)$ to $\mathcal{B}_{k+1}^{(\bullet,\bullet]}$.
A summand $E(a,a;k)$ contributes nothing to either side.
This proves the decomposition of Theorem \ref{thm:main}, the classification of the occurring rectangles into the types (S), (B) and (N) with bottom-right corners on $\Delta$, and the multiset identity of Theorem \ref{thm:dictionary}.
The multiset $\mathcal{B}_k^{(\bullet,\bullet]}(C_*)$ is independent of all choices by the Azumaya--Krull--Remak--Schmidt theorem \cite{Az50}.
Finally, the one-parameter barcodes are recovered from the two-parameter ones by
    \begin{align*}
        \mathcal{B}_k^{\leq\bullet}(C_*)
        ={} & \Set*{[c,\infty) \given R(c;\infty,\infty)\in\mathcal{B}_k^{(\bullet,\bullet]}(C_*)} \\
        & \sqcup \Set*{[c,c+\ell_2) \given R(c;\infty,\ell_2)\in\mathcal{B}_k^{(\bullet,\bullet]}(C_*),\ \ell_2<\infty},
    \end{align*}
    so the two graded barcodes determine each other.
\end{proof}

\begin{remark}[Why the recovery holds]\label{rem:cospan}
    Usher \cite{Ush25} introduced \textit{filtered cospans}, consisting of an ascending chain complex $\mathcal{C}_{\uparrow}$, a descending chain complex $\mathcal{C}^{\downarrow}$, a further chain complex $D_*$, and chain maps $\mathcal{C}_{\uparrow}\to D_*\leftarrow\mathcal{C}^{\downarrow}$, and showed that these stand to interlevel persistence as filtered chain complexes stand to sublevel persistence.
    Every admissible filtered cospan decomposes into eight types of standard elementary summands \cite[Corollary 3.23]{Ush25}, and the ascending elementary complexes used there \cite[Definition 3.16]{Ush25} are precisely the complexes $E(a,b;k)$ and $E(c;k)$ of the present section.
    A Floer-type complex $(C_*,\partial,\ell)$ over $\Lambda^{\FF,\{0\}}$ determines the degenerate filtered cospan in which $\mathcal{C}_{\uparrow}=C_*$ while $\mathcal{C}^{\downarrow}$ and $D_*$ both vanish.
    For such a cospan the decomposition of \cite[Proposition 3.19]{Ush25} leaves only two of the eight types, namely those denoted $(\uparrow_a^b)_k$ and $(\uparrow_a^{\infty})_k$, which are built from $E(a,b;k)$ and $E(c;k)$ respectively.
By Theorem \ref{thm:dictionary} the former accounts for a pair of rectangles of types (B) and (N) in adjacent degrees, and the latter for a rectangle of type (S).
    In particular the summand denoted $(>_a^b)_k$ in \cite{Ush25}, which corresponds to the extended subdiagram of extended persistence and is the only one whose two filtration parameters come from the ascending and the descending complex respectively, does not occur.
    This gives a structural reason for Theorem \ref{thm:dictionary}: genuinely two-parameter information arises from the interaction between an ascending and a descending filtration, whereas the modules $\HH_k^{(\bullet,\bullet]}(C_*)$ are built from a single ascending complex, so that every summand occurring here carries at most one finite filtration parameter.
    Note also that the two-parameter modules of \cite[Section 8]{Ush23} and \cite[Section 6]{Ush25} take their two parameters from the descending and the ascending complex respectively, whereas both of our parameters come from the same ascending complex.
Correspondingly those modules are block-decomposable, while ours is rectangle-decomposable but not block-decomposable.
\end{remark}

\begin{remark}[Relation to homology pyramids]\label{rem:pyramid}
    For a continuous function $f$ on a topological space, the levelset homologies $H_i\bigl(f^{-1}(x,y)\bigr)$ and the relative homologies of the pairs
    \[
        \bigl((x,\infty),(y,\infty)\bigr),\qquad \bigl((-\infty,y),(-\infty,x)\bigr),\qquad \bigl(\RR,(-\infty,x)\cup(y,\infty)\bigr)
    \]
    assemble, after rescaling, into a persistence module over a pyramidal poset, the \textit{homology pyramid} of \cite{CdSM09}.
    It is strongly exact by Mayer--Vietoris and trivial on the two boundary lines of the pyramid, hence interval-decomposable with summands supported on maximal rectangles of the pyramid \cite[Theorem 8.2]{BLO23}.
    The modules $\HH_k^{(\bullet,\bullet]}$ studied here correspond to a single face of that pyramid, namely the one indexed by the pairs $\bigl((-\infty,y),(-\infty,x)\bigr)$ of sublevel sets.
    Consequently, for a Morse function on a closed manifold, Theorem \ref{thm:main} and the classification of the occurring rectangles can be deduced by restricting the decomposition of \cite[Theorem 8.2]{BLO23} to that face.
    Indeed, the relevant maximal rectangles $I_1\times I_2$ are those contained in the open strip $\Set{(x,y) \given |y-x|<1}$, and containment forces $\sup I_2\leq\inf I_1+1$ and $\inf I_2\geq\sup I_1-1$, so that $I_1$ and $I_2$ are bounded and maximality turns both inequalities into equalities.
    The top-left corner $(\inf I_1,\sup I_2)$ therefore lies on the line $y=x+1$ and the bottom-right corner $(\sup I_1,\inf I_2)$ on the line $y=x-1$.
    The latter line is the degenerate edge of the face at hand, on which the pyramid vanishes, and it is carried to $\Delta$ by the identification, which gives the position of the bottom-right corner asserted in Theorem \ref{thm:main}.
    The former line is disjoint from that face, so the restricted rectangle must reach one of the two edges along which the face is glued to the remaining faces, namely the edges where the sublevel set is empty or is everything.
In the coordinates of Theorem \ref{thm:main} these are the limits $\ell_1=\infty$ and $\ell_2=\infty$.
    Hence exactly the types (S), (B) and (N) occur, and the type (R) is excluded.
    The content of Theorem \ref{thm:main} is that the same conclusion holds for an arbitrary Floer-type complex over $\Lambda^{\FF,\{0\}}$, where the remaining faces of the pyramid have no evident analogue.
This is precisely the difficulty that the abstract frameworks of \cite{Ush23,Ush25} are designed to overcome.
    It also explains, from a second point of view, why Theorem \ref{thm:dictionary} holds: a single face of the pyramid cannot carry the extended persistence information that the full pyramid does.
    We thank V.~Lebovici for pointing out this relation.
\end{remark}

\section{The isometry theorem}\label{sec:stability}

In this section we prove Theorem \ref{thm:isometry} and Corollary \ref{cor:distances}.
We first record when the rectangles of Theorem \ref{thm:dictionary} are $2\delta$-trivial.
A half-open rectangle $R(c;\ell_1,\ell_2)=[c-\ell_1,c)\times[c,c+\ell_2)$ as in Theorem \ref{thm:dictionary} is $2\delta$-trivial if and only if $\min\{\ell_1,\ell_2\}\leq 2\delta$.
Hence a rectangle of type (S) is never $2\delta$-trivial for $\delta<\infty$.
A rectangle of type (B) or (N) associated with a finite bar $[a,b)$ via Theorem \ref{thm:dictionary} is $2\delta$-trivial if and only if the bar $[a,b)$ is, that is, if and only if $b-a\leq 2\delta$.

\begin{lemma}\label{lem:direction}
    Let $R=I_1\times I_2$ and $R'=I_1'\times I_2'$ be rectangles in $\RR^2$ such that $\FF_R$ and $\FF_{R'}$ are $\delta$-interleaved and $R$ is not $2\delta$-trivial.
    \begin{enumerate}
        \item If $I_1$ is unbounded below, then $I_1'$ is unbounded below.
        \item If $I_2$ is unbounded above, then $I_2'$ is unbounded above.
    \end{enumerate}
    The analogous statement holds verbatim for interval modules over $\RR$: if $\FF_I$ and $\FF_J$ are $\delta$-interleaved, $I$ is not $2\delta$-trivial, and $I$ is unbounded above or below, then $J$ is unbounded in the same direction.
\end{lemma}

\begin{proof}
Let $F\colon\FF_R\to\FF_{R'}[\vec{\delta}]$ and $G\colon\FF_{R'}\to\FF_R[\vec{\delta}]$ be interleaving morphisms.
We prove (i). The proofs of (ii) and of the one-parameter statement are identical.
Since $R$ is not $2\delta$-trivial, there exists a point $q=(x_0,y_0)\in R$ such that $q+2\vec{\delta}\in R$.
In particular, $y_0,y_0+2\delta\in I_2$.
Since $I_1$ is an interval unbounded below, for every sufficiently negative $x$ we have $x,x+2\delta\in I_1$, so that $p=(x,y_0)$ satisfies $p,p+2\vec{\delta}\in R$.
Then $\bigl(G[\vec{\delta}]\circ F\bigr)_p=\bigl(\Phi_{\FF_R}^{2\vec{\delta}}\bigr)_p=\mathrm{id}_{\FF}\neq 0$, so $F_p\neq 0$ and hence $p+\vec{\delta}\in R'$, and in particular $x+\delta\in I_1'$.
As $x$ can be taken arbitrarily negative, $I_1'$ is unbounded below.
\end{proof}

\begin{corollary}\label{cor:type_preservation}
    Let $R$ and $R'$ be rectangles of types \textup{(S)}, \textup{(B)}, or \textup{(N)} such that $\FF_R$ and $\FF_{R'}$ are $\delta$-interleaved with $\delta<\infty$.
    If the types of $R$ and $R'$ are distinct, then both $R$ and $R'$ are $2\delta$-trivial, and in particular neither is of type \textup{(S)}.
    Similarly, an interval module of a finite bar $[a,b)$ and that of an infinite bar $[c,\infty)$ are never $\delta$-interleaved for $\delta<\infty$.
\end{corollary}

\begin{proof}
Suppose $R$ is not $2\delta$-trivial.
If $R$ is of type (S), Lemma \ref{lem:direction} forces $I_1'$ to be unbounded below and $I_2'$ to be unbounded above, so $R'$ is also of type (S).
If $R$ is of type (B), Lemma \ref{lem:direction} forces $I_1'$ to be unbounded below, so $R'$ is of type (B) or (S). If $R'$ were of type (S), then $R'$ would not be $2\delta$-trivial, and applying Lemma \ref{lem:direction} with the roles of $R$ and $R'$ reversed would force $I_2$ to be unbounded above, contradicting that $R$ is of type (B). Hence $R'$ is of type (B).
The argument for type (N) is symmetric.
Thus a non-$2\delta$-trivial rectangle can only be $\delta$-interleaved with a rectangle of the same type. Therefore, if the types are distinct, $R$ is $2\delta$-trivial, and by symmetry so is $R'$.
A rectangle of type (S) is never $2\delta$-trivial, so it can only be $\delta$-interleaved with another rectangle of type (S).
For the one-parameter statement, $[c,\infty)$ is never $2\delta$-trivial and is unbounded above, so by Lemma \ref{lem:direction} any $\delta$-interleaved partner is unbounded above as well.
\end{proof}

Compare \cite[Proposition 8.5]{Ush25}, which establishes the analogous separation phenomenon, at infinite interleaving distance, for the eight types of standard elementary summands of filtered cospans.

\begin{lemma}\label{lem:interleaving_criterion}
    Let $\delta\in[0,\infty)$.
    \begin{enumerate}
        \item Let $c,d\in\RR$. Then the following are equivalent:
        $\FF_{[c,\infty)}$ and $\FF_{[d,\infty)}$ are $\delta$-interleaved;
        $\FF_{R(c;\infty,\infty)}$ and $\FF_{R(d;\infty,\infty)}$ are $\delta$-interleaved;
        $|c-d|\leq\delta$.
        \item Let $[a,b)$ and $[c,d)$ be finite bars. Then the following are equivalent:
        $\FF_{[a,b)}$ and $\FF_{[c,d)}$ are $\delta$-interleaved;
        $\FF_{R(a;\infty,b-a)}$ and $\FF_{R(c;\infty,d-c)}$ are $\delta$-interleaved;
        $\FF_{R(b;b-a,\infty)}$ and $\FF_{R(d;d-c,\infty)}$ are $\delta$-interleaved;
        $\max\{|a-c|,|b-d|\}\leq\delta$ or both $b-a\leq 2\delta$ and $d-c\leq 2\delta$.
    \end{enumerate}
\end{lemma}

\begin{proof}
We prove that the $\delta$-interleaving of the type (B) rectangle modules is equivalent to the last condition of (ii).
The arguments for the type (N) rectangles (interchange the roles of the two coordinates), for the type (S) rectangles and for the bars themselves (suppress the free coordinate) are identical, and the criterion for bars is classical (see, e.g., \cite[Section 2]{PRSZ20}).

Write $R=(-\infty,a)\times[a,b)$ and $R'=(-\infty,c)\times[c,d)$.

Assume first that $b-a\leq 2\delta$ and $d-c\leq 2\delta$.
Then both rectangles are $2\delta$-trivial and the zero morphisms give a $\delta$-interleaving by Remark \ref{rem:trivial_interleaving}.

Assume next that $\max\{|a-c|,|b-d|\}\leq\delta$.
Define $F\colon\FF_R\to\FF_{R'}[\vec{\delta}]$ by $F_p=\mathrm{id}_{\FF}$ if $p\in R$ and $p+\vec{\delta}\in R'$, and $F_p=0$ otherwise, and define $G$ symmetrically.
We claim that $F$ is a morphism of persistence modules.
Naturality can only fail on a pair $p\preceq q$ for which exactly one of $F_p$, $F_q$ is non-zero while the relevant structure maps are non-zero.
Writing $R''=R'-\vec{\delta}=(-\infty,c-\delta)\times[c-\delta,d-\delta)$, this happens only if either $(\alpha)$ $p\in R\setminus R''$, $q\in R\cap R''$, or $(\beta)$ $p\in R\cap R''$, $q\in R''\setminus R$.
In case $(\alpha)$, we have $p_1\leq q_1<c-\delta$, so the failure $p\notin R''$ forces $p_2\notin[c-\delta,d-\delta)$.
As $p_2\geq a\geq c-\delta$, this means $p_2\geq d-\delta$, whence $q_2\geq p_2\geq d-\delta$ and $q\notin R''$, a contradiction.
In case $(\beta)$, the failure $q\notin R$ forces $q_2\notin[a,b)$.
As $q_2\geq p_2\geq a$ this means $q_2\geq b$, contradicting $q_2<d-\delta\leq b$.
Hence $F$ (and likewise $G$) is a morphism.
For the interleaving relations, let $p,p+2\vec{\delta}\in R$.
Then $p_1+2\delta<a$ gives $p_1+\delta<a-\delta\leq c$, and $p_2\geq a$, $p_2+2\delta<b$ give $p_2+\delta\geq a+\delta\geq c$ and $p_2+\delta<b-\delta\leq d$, so $p+\vec{\delta}\in R'$ and $F_p\neq 0$.
Similarly $G_{p+\vec{\delta}}\neq 0$, and the composition equals $\bigl(\Phi_{\FF_R}^{2\vec{\delta}}\bigr)_p$.
If instead $\Phi^{2\vec{\delta}}$ vanishes at $p$, then so does the composition, since a non-zero composition at $p$ requires $p,p+2\vec{\delta}\in R$.
The second relation is verified symmetrically.

Conversely, assume that $\FF_R$ and $\FF_{R'}$ are $\delta$-interleaved and that they are not both $2\delta$-trivial.
Without loss of generality $b-a>2\delta$.
Evaluating as in the proof of Lemma \ref{lem:direction} at the points $p=(x,y)$ with $x$ arbitrarily negative and $y,y+2\delta\in[a,b)$, we obtain $y+\delta\in[c,d)$ for all $y\in[a,b-2\delta)$, that is,
\begin{equation}\label{eq:onesided}
    [a+\delta,b-\delta)\subset[c,d),\qquad\text{whence}\quad c\leq a+\delta\ \text{and}\ d\geq b-\delta.
\end{equation}
If $d-c>2\delta$, the symmetric argument gives $[c+\delta,d-\delta)\subset[a,b)$, whence $a\leq c+\delta$ and $b\geq d-\delta$.
Combined with \eqref{eq:onesided} this yields $\max\{|a-c|,|b-d|\}\leq\delta$.
If $d-c\leq 2\delta$, then \eqref{eq:onesided} still forces $\max\{|a-c|,|b-d|\}\leq\delta$.
Indeed, if $c<a-\delta$, then
$d\leq c+2\delta<a+\delta<b-\delta$, contradicting $d\geq b-\delta$ from \eqref{eq:onesided}.
Thus $a-\delta\leq c\leq a+\delta$.
Moreover,
$b-\delta\leq d\leq c+2\delta\leq a+3\delta<b+\delta$, where the last inequality again uses $b-a>2\delta$.
Hence $|b-d|\leq\delta$ as well.
\end{proof}

\begin{proof}[Proof of Theorem \textup{\ref{thm:isometry}}]
Fix $\varepsilon\geq 0$ and abbreviate $\mathcal{B}_k^1=\mathcal{B}_k^{\leq\bullet}$ and $\mathcal{B}_k^2=\mathcal{B}_k^{(\bullet,\bullet]}$.
It suffices to prove that, for every $\varepsilon$, there exist $\varepsilon$-matchings between $\mathcal{B}_k^1(C_*)$ and $\mathcal{B}_k^1(D_*)$ for all $k$ if and only if there exist $\varepsilon$-matchings between $\mathcal{B}_k^2(C_*)$ and $\mathcal{B}_k^2(D_*)$ for all $k$.

($\Rightarrow$)
Let $\mu_k\colon\mathcal{B}_k^1(C_*)\nrightarrow\mathcal{B}_k^1(D_*)$ be $\varepsilon$-matchings.
By Corollary \ref{cor:type_preservation}, each $\mu_k$ matches infinite bars with infinite bars and finite bars with finite bars.
Using the bijections of Theorem \ref{thm:dictionary}, define $\nu_k\colon\mathcal{B}_k^2(C_*)\nrightarrow\mathcal{B}_k^2(D_*)$ by transporting the infinite-bar part and the finite-bar part of $\mu_k$ to the rectangles of types (S) and (B), respectively, and the finite-bar part of $\mu_{k-1}$ to the rectangles of type (N).
Matched pairs are $\varepsilon$-interleaved by Lemma \ref{lem:interleaving_criterion}.
Moreover, every rectangle of $\mathcal{B}_k^2$ which is not $2\varepsilon$-trivial corresponds to a bar which is not $2\varepsilon$-trivial, hence lies in the coimage or image, respectively, of the relevant $\mu_j$.
Therefore $\nu_k$ is an $\varepsilon$-matching.

($\Leftarrow$)
Let $\nu_k\colon\mathcal{B}_k^2(C_*)\nrightarrow\mathcal{B}_k^2(D_*)$ be $\varepsilon$-matchings.
By Corollary \ref{cor:type_preservation}, every matched pair consists either of rectangles of the same type, or of two $2\varepsilon$-trivial rectangles.
Discarding the pairs of the latter kind with distinct types does not affect the defining properties of an $\varepsilon$-matching, since $2\varepsilon$-trivial rectangles need not be matched.
Thus we may assume that $\nu_k$ matches only rectangles of the same type.
Define $\mu_k\colon\mathcal{B}_k^1(C_*)\nrightarrow\mathcal{B}_k^1(D_*)$ by transporting the type (S) part of $\nu_k$ to the infinite bars and the type (B) part of $\nu_k$ to the finite bars, via Theorem \ref{thm:dictionary}.
Matched pairs are $\varepsilon$-interleaved by Lemma \ref{lem:interleaving_criterion}.
Every bar of $\mathcal{B}_k^1$ which is not $2\varepsilon$-trivial has a non-$2\varepsilon$-trivial rectangle of type (S) or (B) in $\mathcal{B}_k^2$, which is matched by $\nu_k$ to a rectangle of the same type.
Hence the bar lies in the coimage or image, respectively, of $\mu_k$, and $\mu_k$ is an $\varepsilon$-matching.
\end{proof}

\begin{lemma}\label{lem:slice}
    Let $(C_*,\partial_C,\ell_C)$ and $(D_*,\partial_D,\ell_D)$ be Floer-type complexes over $\Lambda^{\FF,\{0\}}$.
    Then for every $k\in\ZZ$,
    \[
        d_{\mathrm{int}}\left(\HH_k^{\leq\bullet}(C_*),\HH_k^{\leq\bullet}(D_*)\right)
        \leq d_{\mathrm{int}}\left(\HH_k^{(\bullet,\bullet]}(C_*),\HH_k^{(\bullet,\bullet]}(D_*)\right).
    \]
\end{lemma}

\begin{proof}
Let $\delta$ be such that the bipersistence modules are $\delta$-interleaved via $(F,G)$.
Since $C_j$ and $D_j$ are finite-dimensional and orthogonalizable, there is $m\in\RR$ such that $C_j^{\leq t}=D_j^{\leq t}=0$ for all $t<m$ and $j\in\{k-1,k,k+1\}$.
Choose $a_0<m-2\delta$.
Then, for every $b\in\RR$, the quotient complex $C_*^{(a_0,b]}$ coincides with $C_*^{\leq b}$ in the degrees $k-1$, $k$, $k+1$, so that $H_k^{(a_0,b]}(C_*)=H_k^{\leq b}(C_*)$, compatibly with the maps $\iota$.
The same holds for $D_*$, and also with $a_0$ replaced by $a_0+\delta$ or $a_0+2\delta$.
Therefore the components $F_{(a_0,b)}$ and $G_{(a_0,b)}$, $b\in\RR$, define morphisms between $\HH_k^{\leq\bullet}(C_*)$ and the $\delta$-shift of $\HH_k^{\leq\bullet}(D_*)$ and vice versa, and the interleaving relations restrict to the relations $G_{(a_0+\delta,b+\delta)}\circ F_{(a_0,b)}=\iota_{b,b+2\delta}$ and symmetrically, because the comparison maps $\pi_{a_0,a_0+2\delta}^{b+2\delta}$ are the identity in the range above.
Hence the one-parameter modules are $\delta$-interleaved.
\end{proof}

\begin{proof}[Proof of Corollary \textup{\ref{cor:distances}}]
Combining Proposition \ref{thm:converse_algebraic_stability} (applied degreewise), Theorem \ref{thm:isometry}, Theorem \ref{thm:one-para_isometry} (applied degreewise) and Lemma \ref{lem:slice}, we obtain
\begin{align*}
    \sup_k d_{\mathrm{int}}\left(\HH_k^{(\bullet,\bullet]}(C_*),\HH_k^{(\bullet,\bullet]}(D_*)\right)
    &\leq d_{\mathrm{bot}}\left(\mathcal{B}_*^{(\bullet,\bullet]}(C_*),\mathcal{B}_*^{(\bullet,\bullet]}(D_*)\right) \\
    &= d_{\mathrm{bot}}\left(\mathcal{B}_*^{\leq\bullet}(C_*),\mathcal{B}_*^{\leq\bullet}(D_*)\right) \\
    &= \sup_k d_{\mathrm{int}}\left(\HH_k^{\leq\bullet}(C_*),\HH_k^{\leq\bullet}(D_*)\right) \\
    &\leq \sup_k d_{\mathrm{int}}\left(\HH_k^{(\bullet,\bullet]}(C_*),\HH_k^{(\bullet,\bullet]}(D_*)\right),
\end{align*}
so all the inequalities are equalities.
\end{proof}

\begin{remark}\label{rem:constants}
    By Theorem \ref{thm:algebraic_stability}, the constant $3$ in $d_{\mathrm{bot}}\leq 3\,d_{\mathrm{int}}$ is optimal among all p.f.d.\ rectangle-decomposable $\RR^2$-indexed modules.
    Corollary \ref{cor:distances} shows that Floer-type bipersistence modules over $\Lambda^{\FF,\{0\}}$ form a class on which the optimal constant improves to $1$, on par with block-decomposable modules, although the rectangles of types (S), (B) and (N) are not blocks.
\end{remark}

We record a simple padding observation that will be used repeatedly in the applications.

\begin{lemma}[Padding asymmetric shifts]\label{lem:asymmetric_padding}
    Let $\mathfrak{M}$ and $\mathfrak{N}$ be $\RR^n$-indexed persistence modules. Suppose that $e_+,e_-\geq0$ and there are morphisms
    \[
        F\colon\mathfrak{M}\to\mathfrak{N}[\vec{e_+}],
        \qquad
        G\colon\mathfrak{N}\to\mathfrak{M}[\vec{e_-}]
    \]
    satisfying
    \[
        G[\vec{e_+}]\circ F=\Phi_{\mathfrak{M}}^{\vec{\sigma}},
        \qquad
        F[\vec{e_-}]\circ G=\Phi_{\mathfrak{N}}^{\vec{\sigma}},
    \]
    where $\sigma=e_++e_-$.
    Then $\mathfrak{M}$ and $\mathfrak{N}$ are $\eta$-interleaved for $\eta=\max\{e_+,e_-\}$.
\end{lemma}

\begin{proof}
Pad $F$ and $G$ by the structure morphisms from the shifts $e_+$ and $e_-$ to the common shift $\eta$. Naturality of $F$ and $G$ with respect to the structure morphisms, together with the two displayed identities, shows that the two padded compositions are $\Phi_{\mathfrak{M}}^{2\vec\eta}$ and $\Phi_{\mathfrak{N}}^{2\vec\eta}$, respectively.
\end{proof}

\section{Applications in Morse theory and in Floer theory}\label{sec:applications}

\subsection{Morse bipersistence modules}\label{sec:applications_Morse}

Let $M$ be a connected closed manifold and $f\colon M\to\RR$ a Morse function. Fix a Riemannian metric $\rho$ such that $(f,\rho)$ is Morse--Smale, and let $\Cr_k(f)$ denote the set of critical points of $f$ of Morse index $k$.
The Morse complex $(\CM_*(f),\partial^f)$ associated with this Morse--Smale pair is generated over $\Ztwo$ by the critical points of $f$, graded by the Morse index, with $\partial^f$ counting negative gradient flow lines.
Together with the filtration function $\ell\bigl(\sum_i n_ip_i\bigr)=\max\{\,f(p_i)\mid n_i\neq 0\,\}$, it is a Floer-type complex over $\Lambda^{\Ztwo,\{0\}}=\Ztwo$ \cite[Example 2.5]{UZ16} (see \cite{Sch93,BH04,AD14} for details on Morse homology). Changing the Morse--Smale metric gives filtered continuation isomorphisms with zero filtration shift; alternatively, the filtered and action-window Morse homologies can be identified with the singular homologies of sublevel sets and of pairs of sublevel sets, respectively. Hence the persistence modules and barcodes below are independent of the choice of $\rho$.

\begin{definition}
    Given a Morse function $f$ on $M$ and $k\in\ZZ$, the p.f.d.\ Floer-type bipersistence module $\HH_k^{(\bullet,\bullet]}\bigl(\CM_*(f)\bigr)$ is called the \textit{Morse bipersistence module in degree $k$}, denoted by $\HHMM_k^{(\bullet,\bullet]}(f)$.
\end{definition}

By Theorem \ref{thm:main}, Morse bipersistence modules are rectangle-decomposable, and we write $\mathcal{B}_k^{(\bullet,\bullet]}(f)$ and $\mathcal{B}_*^{(\bullet,\bullet]}(f)$ for the associated rectangle barcode and its graded version.
Similarly, let $\HHMM_k^{\leq\bullet}(f)$ and $\mathcal{B}_k^{\leq\bullet}(f)$ denote the one-parameter Morse persistence module and its barcode.

The next theorem shows that, at each critical value and in each degree, the number of rectangle summands equals the number of critical points.

\begin{theorem}[Morse Generator--Rectangle Theorem]\label{thm:morse_normal_form}
    Let $M$ be a connected closed manifold and $f\colon M\to\RR$ a Morse function.
    For each degree $k\in\ZZ$, there exists a family of rectangles $R_p=R(f(p);\ell_{1p},\ell_{2p})$, indexed by $p\in\Cr_k(f)$, such that
    \[
        \HHMM_k^{(\bullet,\bullet]}(f) \cong \bigoplus_{p\in\Cr_k(f)} (\Ztwo)_{R_p}.
    \]
    Every $R_p$ is of type $(\mathrm{S})$, $(\mathrm{B})$, or $(\mathrm{N})$, and the multiset $\{R_p\mid p\in\Cr_k(f)\}$ is uniquely determined by $f$. If several critical points of index $k$ have the same critical value, the indexing of the corresponding rectangles by those critical points need not be unique, although the indexing is unique whenever the critical values in degree $k$ are pairwise distinct.
\end{theorem}

\begin{proof}
First, every corner coordinate is a critical value in the relevant degree. Indeed, if $c$ is not a critical value, then for all sufficiently small $\varepsilon>0$ the quotient complex $\CM_*^{(c-\varepsilon,c]}(f)$ is zero. Hence $\HM_k^{(c-\varepsilon,c]}(f)=0$, whereas a rectangle with corner $(c,c)$ would contribute non-trivially to this group. Thus a rectangle in degree $k$ with corner $(c,c)$ forces the presence of an index-$k$ critical point at level $c$.
Now fix a critical value $c$ and choose $\varepsilon>0$ so small that $(c-\varepsilon,c)$ contains no critical value. The action-window Morse homology decomposes into local Morse homology groups (see, e.g., \cite[Section 3.1]{GG10}):
\[
    \HM_k^{(c-\varepsilon,c]}(f) \cong \bigoplus_{\substack{p\in\Cr_k(f)\\ f(p)=c}} \HM_k^{\mathrm{loc}}(p).
\]
Since each non-degenerate critical point has one-dimensional local Morse homology in its Morse degree, the dimension of this group equals the number of critical points $p\in\Cr_k(f)$ with $f(p)=c$. On the other hand, by Theorem \ref{thm:dictionary} a finite horizontal side length of a rectangle with corner $(c,c)$ has the form $c-a$ for a bar $[a,c)$ of $\mathcal{B}_{k-1}^{\leq\bullet}(f)$, so that $a$ is again a critical value. Since $(c-\varepsilon,c)$ contains no critical value, every rectangle with corner $(c,c)$ satisfies $\ell_1\geq\varepsilon$ and hence contains the point $(c-\varepsilon,c)$. Conversely, if a rectangle with corner $(d,d)$ contains this point, then $c-\varepsilon<d\leq c$; the first part of the proof and the choice of $\varepsilon$ force $d=c$. Therefore the rectangle decomposition shows that the same dimension is precisely the multiplicity of rectangles in $\mathcal{B}_k^{(\bullet,\bullet]}(f)$ with bottom-right corner $(c,c)$. Hence these two multiplicities agree for every $c$. Choosing, for each $c$, a bijection between the critical points at level $c$ and the rectangles with corner $(c,c)$ gives the asserted indexing. The multiset of rectangles is the barcode and is therefore unique, whereas the chosen bijection need not be unique when several critical points have the same critical value.
\end{proof}

\begin{remark}[Barannikov pairing]\label{rem:barannikov}
    Combining Theorems \ref{thm:morse_normal_form} and \ref{thm:dictionary}, the multiset of finite side lengths attached to a fixed critical value is explicitly determined by the one-parameter barcodes. If the critical values in each degree are pairwise distinct, every critical point $p$ of index $k$ is uniquely of one of the following kinds: it opens a bar of $\mathcal{B}_k^{\leq\bullet}(f)$ at $f(p)$, giving type (S) when the bar is infinite and type (B) with $\ell_{2p}$ equal to the bar length otherwise, or it closes a bar $[a,f(p))$ of $\mathcal{B}_{k-1}^{\leq\bullet}(f)$, giving type (N) with $\ell_{1p}=f(p)-a$. When several critical points share the same critical value and degree, the same description holds after choosing an indexing as in Theorem \ref{thm:morse_normal_form}, but the assignment to the individual critical points need not be canonical.
This is the barcode form of Barannikov's classical pairing \cite{barannikov1994}.
    In the one-parameter barcode, the generators corresponding to type (N) rectangles are visible through right endpoints of bars one degree below, whereas the rectangle barcode records the same multiplicity directly in the higher degree.
\end{remark}

We conclude the Morse-theoretic discussion with functional stability.

\begin{proof}[Proof of Theorem \textup{\ref{thm:functional_stability}}]
Put $\delta=\|f-g\|$.
Writing $M_f^t=f^{-1}((-\infty,t])$ and $M_g^t=g^{-1}((-\infty,t])$, for every $t\in\RR$ the inclusions
\[
    M_f^t\subset M_g^{t+\delta},
    \qquad
    M_g^t\subset M_f^{t+\delta}
\]
induce a $\delta$-interleaving between the one-parameter Morse persistence modules in every degree.
The one-parameter isometry theorem (Theorem \ref{thm:one-para_isometry}), followed by Theorem \ref{thm:isometry}, therefore gives
\[
    d_{\mathrm{bot}}\left(\mathcal{B}_*^{(\bullet,\bullet]}(f),\mathcal{B}_*^{(\bullet,\bullet]}(g)\right)
    =d_{\mathrm{bot}}\left(\mathcal{B}_*^{\leq\bullet}(f),\mathcal{B}_*^{\leq\bullet}(g)\right)
    \leq\delta.
\]
To see that the constant is sharp, take $g=f+c$ with $c>0$.
Since $M$ is connected, the degree-zero barcode of $f$ has a unique infinite bar, namely $[\min_M f,\infty)$, while that of $g$ has the unique infinite bar $[\min_M f+c,\infty)$.
These bars must be matched, at cost $c$, in every finite bottleneck matching.
The reverse inequality follows, and hence Theorem \ref{thm:isometry} gives
\[
    d_{\mathrm{bot}}\left(\mathcal{B}_*^{(\bullet,\bullet]}(f),\mathcal{B}_*^{(\bullet,\bullet]}(g)\right)
    =c=\|f-g\|.
\]
\end{proof}

\subsection{Hamiltonian bipersistence modules: contractible orbits}\label{sec:filtered_Floer_contractible}

Let $(M,\omega)$ be a connected closed symplectic manifold.
Let $\LL_\alpha M$ denote the component of the free loop space of $M$ corresponding to a free homotopy class $\alpha\in[S^1,M]$, with a fixed reference loop $z_\alpha$, and let $0$ denote the class of contractible loops.
Evaluating $[\omega]$ and the first Chern class $c_1(TM)$ on spheres defines homomorphisms $I_\omega,I_{c_1}\colon\pi_2(M)\to\RR$, and evaluating them on tori defines homomorphisms $I_\omega^\alpha,I_{c_1}^\alpha\colon\pi_1(\LL_\alpha M,z_\alpha)\to\RR$.

\begin{definition}\label{def:atoroidal_aspherical}
    A connected closed symplectic manifold $(M,\omega)$ is \textit{symplectically aspherical} if $I_\omega=I_{c_1}=0$, and \textit{symplectically atoroidal} if $I_\omega^\alpha=I_{c_1}^\alpha=0$ for all $\alpha\in[S^1,M]$.
\end{definition}

Symplectically atoroidal manifolds are symplectically aspherical.

A Hamiltonian $H\colon S^1\times M\to\RR$, $H_t=H(t,\cdot)$, generates the Hamiltonian isotopy $\{\varphi_H^t\}$ of the vector field $X_{H_t}$ defined by $\iota_{X_{H_t}}\omega=-dH_t$, with time-one map $\varphi_H$.
The group of Hamiltonian diffeomorphisms is denoted by $\Ham(M)$.
Let $\mathcal{P}_1(H;\alpha)$ denote the set of one-periodic orbits of the isotopy in the class $\alpha$.
An orbit $x\in\mathcal{P}_1(H;\alpha)$ is \textit{non-degenerate} if $1$ is not an eigenvalue of $(d\varphi_H)_{x(0)}$, and $H$ is \textit{non-degenerate} if all its one-periodic orbits are.
We call $H$ \textit{normalized} if $\int_M H_t\,\omega^n=0$ for all $t$.
Given a Hamiltonian $H$, put
\begin{equation}\label{eq:E^+_E^-}
    E^+(H)=\int_0^1\max_M H_t\,dt,\qquad E^-(H)=-\int_0^1\min_M H_t\,dt,
\end{equation}
and $\osc H = E^+(H)+E^-(H)$.
The \textit{Hofer norm} of $\phi\in\Ham(M)$ is $\|\phi\|_{\mathrm{Hofer}}=\inf\Set{\osc H \given \varphi_H=\phi}$, and the \textit{Hofer metric} is $d_{\mathrm{Hofer}}(\phi,\psi)=\|\phi^{-1}\circ\psi\|_{\mathrm{Hofer}}$, a bi-invariant genuine metric on $\Ham(M)$ \cite{Ho90,Po93,LM95}.

Assume from now on that $(M,\omega)$ is symplectically aspherical.
Then the \textit{action} of a contractible loop $x$ with a capping disc $v$,
\[
    \mathcal{A}_H(x)=-\int_{D^2}v^*\omega+\int_0^1 H_t\bigl(x(t)\bigr)\,dt,
\]
does not depend on the capping, and every non-degenerate orbit $x\in\mathcal{P}_1(H;0)$ carries a well-defined Conley--Zehnder index $\mu_{\mathrm{CZ}}(H,x)\in\ZZ$ (see, e.g., \cite{SZ92}).
For a non-degenerate $H$ and a generic loop $J$ of compatible almost complex structures, the Floer complex $(\CF_*(H),\partial^{H,J})$ is generated over $\Ztwo$ by $\mathcal{P}_1(H;0)$, graded by $\mu_{\mathrm{CZ}}$, with $\partial^{H,J}$ counting Floer trajectories of relative index one.
The filtration function $\ell(\xi)=\max\{\,\mathcal{A}_H(x)\mid\xi_{x}\neq 0\,\}$ makes $(\CF_*(H),\partial^{H,J},\ell)$ a Floer-type complex over $\Lambda^{\Ztwo,\{0\}}=\Ztwo$ \cite[Example 4.2]{UZ16} (see \cite{AD14,MS12} for detailed expositions). Changing the generic loop $J$ gives filtered continuation isomorphisms with zero Hamiltonian change, so the persistence modules below are independent of $J$ and we suppress it from the notation. Moreover, the action-preserving naturality isomorphisms of Schwarz \cite{Sch00} show that, for normalized Hamiltonians on a symplectically aspherical manifold, the resulting filtered Floer homology depends only on the time-one map $\phi=\varphi_H$.

\begin{definition}
    Let $(M,\omega)$ be symplectically aspherical and $\phi=\varphi_H\in\Ham(M)$ a non-degenerate Hamiltonian diffeomorphism generated by a normalized Hamiltonian $H$.
    Given a degree $k$, the Floer-type bipersistence module $\HH_k^{(\bullet,\bullet]}\bigl(\CF_*(H)\bigr)$ is called the \textit{Hamiltonian bipersistence module in degree $k$}, denoted by $\HHFF_k^{(\bullet,\bullet]}(\phi)$.
\end{definition}

By Theorem \ref{thm:main}, the rectangle barcode $\mathcal{B}_k^{(\bullet,\bullet]}(\phi)$ is defined, and we write $\mathcal{B}_*^{(\bullet,\bullet]}(\phi)$ for its graded version.

\begin{theorem}[Floer Orbit--Rectangle Theorem: contractible orbits]\label{thm:floer_contractible_normal_form}
    Let $(M,\omega)$ be a symplectically aspherical symplectic manifold and $\phi=\varphi_H$ a non-degenerate Hamiltonian diffeomorphism generated by a normalized Hamiltonian $H$.
    For each degree $k\in\ZZ$, there exists a family of rectangles $R_x=R(\mathcal{A}_H(x);\ell_{1x},\ell_{2x})$, indexed by the orbits $x\in\mathcal{P}_1(H;0)$ with $\mu_{\mathrm{CZ}}(H,x)=k$, such that
    \[
        \HHFF_k^{(\bullet,\bullet]}(\phi) \cong \bigoplus_{\substack{x\in\mathcal{P}_1(H;0)\\ \mu_{\mathrm{CZ}}(H,x)=k}} (\Ztwo)_{R_x}.
    \]
    Every $R_x$ is of type $(\mathrm{S})$, $(\mathrm{B})$, or $(\mathrm{N})$, and the multiset of these rectangles is uniquely determined by $\phi$. If several degree-$k$ orbits have the same action, the indexing by the individual orbits need not be unique, although it is unique whenever the action values in degree $k$ are pairwise distinct.
\end{theorem}

\begin{proof}
First, every corner coordinate belongs to the action spectrum in the relevant degree. If $c$ is not an action value, then, since the action spectrum of a non-degenerate Hamiltonian on a closed manifold is finite, the quotient complex $\CF_*^{(c-\varepsilon,c]}(H)$ is zero for all sufficiently small $\varepsilon>0$. Hence $\HF_k^{(c-\varepsilon,c]}(H)=0$, whereas a rectangle in degree $k$ with corner $(c,c)$ would contribute non-trivially to this group. Thus such a rectangle forces a degree-$k$ one-periodic orbit of action $c$.
Now fix an action value $c$ and choose $\varepsilon>0$ so small that $(c-\varepsilon,c)$ contains no action value. Since $H$ is non-degenerate, the action-window Floer homology decomposes into local Floer homology groups (see \cite[Section 3.2]{GG10}):
\[
    \HF_k^{(c-\varepsilon,c]}(H) \cong \bigoplus_{\substack{x\in\mathcal{P}_1(H;0)\\ \mu_{\mathrm{CZ}}(H,x)=k\\ \mathcal{A}_H(x)=c}} \HF_k^{\mathrm{loc}}(x).
\]
Each non-degenerate orbit has one-dimensional local Floer homology in its Conley--Zehnder degree. Comparing dimensions with the rectangle decomposition, exactly as in the proof of Theorem \ref{thm:morse_normal_form}, shows that the multiplicity of rectangles in degree $k$ with corner $(c,c)$ equals the number of degree-$k$ orbits of action $c$. Choosing a bijection at each action value yields the asserted indexing, with the stated non-uniqueness when action values coincide.
\end{proof}

By the same dictionary as in Remark \ref{rem:barannikov}, the multiset of finite side lengths at each action value is determined by the one-parameter barcodes. When the action values in each degree are pairwise distinct, every one-periodic orbit can therefore be assigned uniquely to the bar it opens or closes, whereas without this simplicity assumption the assignment is only determined up to a bijection among orbits with the same action and degree.

\begin{proof}[Proof of Theorem \textup{\ref{thm:dynamical_stability}}]
Fix $\varepsilon>0$. Choose a normalized Hamiltonian $K$ generating $\phi$ and a normalized Hamiltonian $L$ generating $\phi^{-1}\psi$ such that
\[
    \osc L<d_{\mathrm{Hofer}}(\phi,\psi)+\varepsilon.
\]
Let $H=K\#L$ generate the composed path $\varphi_K^t\circ\varphi_L^t$.
Then $H$ is normalized and $\varphi_H=\psi$, and
\[
    H_t-K_t=L_t\circ(\varphi_K^t)^{-1},
\]
and therefore $E^\pm(H-K)=E^\pm(L)$ and $\osc(H-K)=\osc L$.
The standard continuation maps from $K$ to $H$ and back have shifts $e_+=E^+(H-K)$ and $e_-=E^+(K-H)$, and their compositions are the comparison morphisms of shift $e_++e_-$. Lemma \ref{lem:asymmetric_padding} yields an $\eta$-interleaving, where $\eta=\max\{e_+,e_-\}\leq\osc(H-K)$. Hence Corollary \ref{cor:distances} gives
\[
    d_{\mathrm{bot}}\left(\mathcal{B}_*^{(\bullet,\bullet]}(\phi),\mathcal{B}_*^{(\bullet,\bullet]}(\psi)\right)
    \leq \osc L
    < d_{\mathrm{Hofer}}(\phi,\psi)+\varepsilon.
\]
Letting $\varepsilon\to0$ proves the theorem.
\end{proof}

We next make the classical invariants visible in the rectangle barcode.
Let $\HHFF_k^{\leq\bullet}(\phi)=\HH_k^{\leq\bullet}\bigl(\CF_*(H)\bigr)$ denote the one-parameter Floer persistence module and $\mathcal{B}_*^{\leq\bullet}(\phi)$ its graded barcode, and let $\HF(\phi)$ denote the total homology of $\CF_*(H)$.
The \textit{spectral invariant} $c\colon\HF(\phi)\setminus\{0\}\to\RR$ is defined by
\[
    c(\beta)=\inf\Set*{b\in\RR \given \beta\in\Image\bigl(H_*^{\leq b}\to\HF(\phi)\bigr)},
\]
a characteristic exponent of $\HHFF_*^{\leq\bullet}(\phi)$ in the sense of \cite[Section 4.1.1]{PRSZ20}.
Its image is the set of left endpoints of the infinite bars in $\mathcal{B}_*^{\leq\bullet}(\phi)$ \cite{Sch00,Oh05,Vi92}; for this barcode formulation see also \cite{UZ16,PRSZ20}.
The \textit{boundary depth} of $\phi=\varphi_H$ \cite{Ush11,Ush13} is
\[
    \beta(\phi)=\inf\Set*{b\geq 0 \given (\Image\partial^H)\cap\CF^{\leq t}(H)\subset\partial^H\bigl(\CF^{\leq t+b}(H)\bigr)\ \text{for all}\ t\in\RR},
\]
that is, the maximum of $0$ and the lengths of the finite bars in $\mathcal{B}_*^{\leq\bullet}(\phi)$.

By Theorem \ref{thm:dictionary}, the rectangles of type (S) correspond bijectively to the infinite bars $[c,\infty)$ in $\mathcal{B}_*^{\leq\bullet}(\phi)$, while the rectangles of type (B) correspond bijectively to the finite bars $[c,c+\ell_2)$.
Hence the image of the spectral invariant $c\colon\HF(\phi)\setminus\{0\}\to\RR$ coincides with the set $\Set[\big]{c\in\RR \given R(c;\infty,\infty)\in\mathcal{B}_*^{(\bullet,\bullet]}(\phi)}$ of diagonal coordinates of the corners of the rectangles of type (S), and the boundary depth is the maximum of $0$ and the finite vertical side lengths of the rectangles of type $(\mathrm{B})$:
\[
    \beta(\phi)=\max\left(\{0\}\cup\Set*{\ell_2\in(0,\infty) \given R(c;\infty,\ell_2)\in\mathcal{B}_*^{(\bullet,\bullet]}(\phi)\ \text{for some}\ c}\right).
\]

In analogy with the boundary depth, the earlier version of this paper introduced the \textit{non-cycle depth}
\[
    \alpha(\phi)=\max\left(\{0\}\cup\Set*{\ell_1\in(0,\infty) \given R(c;\ell_1,\infty)\in\mathcal{B}_*^{(\bullet,\bullet]}(\phi)\ \text{for some}\ c}\right),
\]
the maximum of $0$ and the finite horizontal side lengths of the rectangles of type (N), and asked whether it carries new information.
Theorem \ref{thm:dictionary} answers this question.

\begin{corollary}\label{cor:noncycle_depth}
    Let $(M,\omega)$ be symplectically aspherical and $\phi\in\Ham(M)$ non-degenerate.
    Then $\alpha(\phi)=\beta(\phi)$.
\end{corollary}

\begin{proof}
By Theorem \ref{thm:dictionary}, the rectangles of type (N) in $\mathcal{B}_*^{(\bullet,\bullet]}(\phi)$ are exactly the rectangles $R(b;b-a,\infty)$ associated with the finite bars $[a,b)$ of $\mathcal{B}_*^{\leq\bullet}(\phi)$, one degree above.
Hence the multiset of finite horizontal side lengths of the rectangles of type (N) coincides with the multiset of lengths of the finite bars, which in turn coincides with the multiset of finite vertical side lengths of the rectangles of type (B).
Taking maxima yields $\alpha(\phi)=\beta(\phi)$.
\end{proof}

\subsection{Hamiltonian bipersistence modules: non-contractible orbits}\label{sec:filtered_Floer_non-contractible}

Assume now that $(M,\omega)$ is symplectically atoroidal and fix a non-trivial class $\alpha\in[S^1,M]$. Fix the reference loop $z_\alpha$ introduced above and a symplectic trivialization of $TM$ along $z_\alpha$; this reference data is used throughout the subsection.
For a loop $x\in\LL_\alpha M$ and a cylinder $\Pi$ connecting $z_\alpha$ to $x$, the action
$\mathcal{A}_{H,\alpha}(x)=-\int_{[0,1]\times S^1}\Pi^*\omega+\int_0^1 H_t\bigl(x(t)\bigr)\,dt$
does not depend on $\Pi$, and the chosen trivialization gives well-defined Conley--Zehnder indices of non-degenerate orbits.
For a non-degenerate $H$ one obtains, as in the contractible case, the Floer-type complex $(\CF_*(H;\alpha),\partial^{H,J},\ell)$ over $\Ztwo$ and the Floer homology $\HF_*(H;\alpha)$. Changing $J$ gives filtered continuation isomorphisms with zero Hamiltonian change, so the associated persistence modules are independent of this auxiliary choice. See \cite{Gu13,GG16,PS16} for detailed expositions.

\begin{remark}\label{rem:non-contractible_vanishing}
    The total Floer homology $\HF_*(H;\alpha)$ vanishes for $\alpha\neq 0$, since Floer homology is independent of the Hamiltonian and every $C^2$-small autonomous Hamiltonian has only contractible periodic orbits.
    Hence spectral invariants are trivial in this setting, and it is essential to work with action windows.
\end{remark}

\begin{definition}
    Let $(M,\omega)$ be symplectically atoroidal, $H$ a non-degenerate Hamiltonian, and $\alpha\neq 0$.
    Given a degree $k$, we call $\HH_k^{(\bullet,\bullet]}\bigl(\CF_*(H;\alpha)\bigr)$ the \textit{Hamiltonian bipersistence module in the class $\alpha$ in degree $k$} and denote it by $\HHFF_k^{(\bullet,\bullet]}(H;\alpha)$.
    By Theorem \ref{thm:main}, its rectangle barcode $\mathcal{B}_k^{(\bullet,\bullet]}(H;\alpha)$ is defined.
\end{definition}

\begin{theorem}[Floer Orbit--Rectangle Theorem: non-contractible orbits]\label{thm:floer_non-contractible_normal_form}
    Let $(M,\omega)$ be symplectically atoroidal, $\alpha\in[S^1,M]$ non-trivial, and $H$ a non-degenerate Hamiltonian.
    For each degree $k\in\ZZ$, there exists a family of rectangles $R_x=R(\mathcal{A}_{H,\alpha}(x);\ell_{1x},\ell_{2x})$, indexed by the orbits $x\in\mathcal{P}_1(H;\alpha)$ with $\mu_{\mathrm{CZ}}(H,x)=k$, such that
    \[
        \HHFF_k^{(\bullet,\bullet]}(H;\alpha) \cong \bigoplus_{\substack{x\in\mathcal{P}_1(H;\alpha)\\ \mu_{\mathrm{CZ}}(H,x)=k}} (\Ztwo)_{R_x}.
    \]
    Every $R_x$ is of type $(\mathrm{B})$ or $(\mathrm{N})$, and the multiset of these rectangles is uniquely determined. If several degree-$k$ orbits have the same action, the indexing by the individual orbits need not be unique, although it is unique whenever the action values in degree $k$ are pairwise distinct.
\end{theorem}

\begin{proof}
The proof of Theorem \ref{thm:floer_contractible_normal_form} applies verbatim.
By Remark \ref{rem:non-contractible_vanishing} and Theorem \ref{thm:dictionary}, the one-parameter barcode $\mathcal{B}_*^{\leq\bullet}(H;\alpha)$ contains no infinite bar, so no rectangle of type (S) occurs.
\end{proof}

\begin{theorem}[Hamiltonian Stability Theorem]\label{thm:dynamical_stability_non-contractible}
    Let $(M,\omega)$ be symplectically atoroidal and $\alpha\in[S^1,M]$ non-trivial. For any pair of normalized non-degenerate Hamiltonians $K$ and $H$, we have
    \[
        d_{\mathrm{bot}}\left(\mathcal{B}_*^{(\bullet,\bullet]}(K;\alpha),\mathcal{B}_*^{(\bullet,\bullet]}(H;\alpha)\right) \leq \max\bigl\{E^+(H-K),E^+(K-H)\bigr\},
    \]
    where $E^+$ is as in \eqref{eq:E^+_E^-}.
    In particular, the left-hand side is at most $\osc(H-K)$.
\end{theorem}

\begin{proof}
Put $e_+=E^+(H-K)$ and $e_-=E^+(K-H)$. Since $K$ and $H$ are normalized, $e_+,e_-\geq0$. The standard continuation maps in the class $\alpha$ have filtration shifts $e_+$ and $e_-$, and the continuation-composition argument identifies their two compositions with the comparison morphisms of shift $e_++e_-$. Lemma \ref{lem:asymmetric_padding} therefore gives an $\eta$-interleaving in every degree, where
\[
    \eta=\max\{e_+,e_-\}\leq e_++e_-=\osc(H-K).
\]
Corollary \ref{cor:distances} gives the stated bottleneck bound.
\end{proof}

\begin{remark}
The normalization is essential for a barcode estimate of this kind: adding a time-dependent spatial constant to a Hamiltonian leaves its Hamiltonian isotopy unchanged but translates all action values, and hence translates the rectangle barcode along the diagonal. We do not claim here that the non-contractible rectangle barcode depends only on the time-one map, since such a statement would require an additional naturality argument with respect to the Hamiltonian path.
\end{remark}

We now turn to the spectral spread, inspired by Polterovich and Shelukhin's $\ZZ/k\ZZ$ spectral spreads \cite[Section 4.1]{PS16}.
Given $d\geq 0$ and an action window $(a,b]$, let
$\Phi_H^{(a,b],d} = \pi_{a,a+d}^{b+d}\circ\iota_{b,b+d}^a\colon \HF_*^{(a,b]}(H;\alpha)\to \HF_*^{(a+d,b+d]}(H;\alpha)$
denote the comparison map.

\begin{definition}\label{def:Floer_spectral_spread}
    Let $(M,\omega)$ be symplectically atoroidal, $\alpha\in[S^1,M]$ non-trivial, and $H$ a non-degenerate Hamiltonian. The \textit{spectral spread of $H$ in $\alpha$} is
    \[
        w_\alpha(H)=\sup\left(\{0\}\cup\Set*{d\geq 0 \given \Phi_H^{(a,b],d}\neq 0\ \text{for some action window}\ (a,b]}\right).
    \]
\end{definition}

The quantity $w_\alpha(H)$ is finite. Indeed, a non-degenerate Hamiltonian on a closed manifold has only finitely many one-periodic orbits, and Remark \ref{rem:non-contractible_vanishing} together with Theorem \ref{thm:dictionary} shows that its finite barcode contains only rectangles of types (B) and (N). The diagonal structure maps of each such summand vanish once the shift reaches its finite side length. Proposition \ref{prop:Floer_spectral_spread} below gives the precise formula.

\begin{proposition}[{cf.\ \cite[Proposition 4.2]{PS16}}]\label{prop:spectral_spread_properties}
    Let $(M,\omega)$ be symplectically atoroidal, let $\alpha\in[S^1,M]$ be non-trivial, and let $K,H$ be normalized non-degenerate Hamiltonians. Then
    \[
        |w_\alpha(K)-w_\alpha(H)|\leq \osc(H-K).
    \]
\end{proposition}

\begin{proof}
Put $e_+=E^+(H-K)$, $e_-=E^+(K-H)$, and $\delta=e_++e_-=\osc(H-K)$.
Let $\sigma_{KH}$ and $\sigma_{HK}$ denote the continuation morphisms, with filtration shifts $e_+$ and $e_-$, respectively.
For every action window $I$, naturality with respect to the comparison maps and the continuation-composition identity give
\[
    \sigma_{HK}\circ\Phi_H^{I+e_+,d}\circ\sigma_{KH}=\Phi_K^{I,\delta+d},
\]
where $I+e_+$ denotes the translated action window.
Let $r<w_\alpha(K)$.
If $r\leq\delta$, then $r\leq w_\alpha(H)+\delta$ because $w_\alpha(H)\geq0$.
If $r>\delta$, there exist $d'>r$ and an action window $I$ with $\Phi_K^{I,d'}\neq0$; since $\Phi_K^{I,d'}$ factors through $\Phi_K^{I,r}$, we have $\Phi_K^{I,r}\neq0$. Put $d=r-\delta>0$.
The displayed identity then forces $\Phi_H^{I+e_+,d}\neq0$, and hence $r-\delta=d\leq w_\alpha(H)$.
Letting $r\nearrow w_\alpha(K)$ gives $w_\alpha(K)\leq w_\alpha(H)+\delta$.
Interchanging $K$ and $H$ gives the reverse inequality.
\end{proof}

The spectral spread can be read off from the rectangle barcode.

\begin{proposition}\label{prop:Floer_spectral_spread}
    For every non-degenerate Hamiltonian $H$,
    \[
        w_\alpha(H)=\max\left(\{0\}\cup\Set*{\min\{\ell_{1R},\ell_{2R}\} \given R\in\mathcal{B}_*^{(\bullet,\bullet]}(H;\alpha)}\right),
    \]
    where $\ell_{1R},\ell_{2R}\in(0,\infty]$ denote the side lengths of a rectangle $R=R(c;\ell_{1R},\ell_{2R})$.
\end{proposition}

\begin{proof}
Decompose the bipersistence module into rectangle summands. The comparison map $\Phi_H^{(a,b],d}$ is non-zero if and only if some rectangle $R$ in the barcode contains both $(a,b)$ and $(a+d,b+d)$. For $R=R(c;\ell_{1R},\ell_{2R})$, such a pair of points exists exactly when $d<\min\{\ell_{1R},\ell_{2R}\}$. Since $H$ is non-degenerate on a closed manifold, the barcode is finite. Taking the supremum over $d$ for each rectangle and then the maximum over the finitely many rectangles gives the formula, with the adjoined value $0$ covering the case of an empty barcode.
\end{proof}

Since no rectangle of type (S) occurs (Theorem \ref{thm:floer_non-contractible_normal_form}), Proposition \ref{prop:Floer_spectral_spread} and Theorem \ref{thm:dictionary} show that $w_\alpha(H)$ equals the boundary depth of the Floer-type complex $\CF_*(H;\alpha)$ in the sense of \cite{Ush11,UZ16}, that is, $0$ when the barcode is empty and otherwise the largest finite bar length.
This identification was observed by V.~Lebovici and by Q.~Feng (see the acknowledgments).

\begin{proposition}\label{prop:existence_of_non-contractible_orbits}
    Let $(M,\omega)$ be symplectically atoroidal, $\alpha\in[S^1,M]$ non-trivial, and let $K,H$ be normalized non-degenerate Hamiltonians. If
    \[
        \osc(H-K)<w_\alpha(K),
    \]
    then $H$ has a one-periodic orbit in the class $\alpha$.
    Equivalently, by Proposition \ref{prop:Floer_spectral_spread}, it is enough that the barcode $\mathcal{B}_*^{(\bullet,\bullet]}(K;\alpha)$ contain a rectangle $R=R(c;\ell_{1R},\ell_{2R})$ with
    \[
        \osc(H-K)<\min\{\ell_{1R},\ell_{2R}\}.
    \]
\end{proposition}

\begin{proof}
By Proposition \ref{prop:spectral_spread_properties},
\[
    0<w_\alpha(K)-\osc(H-K)\leq w_\alpha(H).
\]
Hence $\Phi_H^{I,d}\neq 0$ for some window $I$ and some $d>0$. In particular $\HF_*^I(H;\alpha)\neq 0$, so $H$ has a one-periodic orbit in the class $\alpha$. The barcode formulation follows from Proposition \ref{prop:Floer_spectral_spread}.
\end{proof}

\begin{remark}[Relaxing the hypotheses]\label{rem:relaxing_hypotheses}
    Two hypotheses above can be relaxed.
    First, the normalization of $K$ and $H$ in Propositions \ref{prop:spectral_spread_properties} and \ref{prop:existence_of_non-contractible_orbits} is inessential: replacing $K$ and $H$ by their normalizations changes neither the generated isotopies and their one-periodic orbits nor $\osc(H-K)$, while it translates all action values of each Hamiltonian by a common constant, so that the rectangle barcodes translate along the diagonal and the spectral spreads are unchanged.
    Second, the non-degeneracy of $H$ in Proposition \ref{prop:existence_of_non-contractible_orbits} can be dropped: since non-degenerate Hamiltonians are dense in the $C^2$-topology, one may choose non-degenerate Hamiltonians $H_i$ converging to $H$ in that topology with $\osc(H_i-K)<w_\alpha(K)$, and each $H_i$ then has a one-periodic orbit $x_i$ in the class $\alpha$.
    Since $X_{H_i}\to X_H$ in $C^1$, the loops $x_i$ are equi-Lipschitz; hence, by the Arzel\`a--Ascoli theorem, after passing to a subsequence they converge uniformly to a loop $x$, and the integral form of Hamilton's equation shows that $x$ is a one-periodic orbit of $H$.
    For all sufficiently large $i$, the loops $x_i$ and $x$ are uniformly closer than the injectivity radius of a fixed Riemannian metric on $M$, hence freely homotopic; therefore $x$ lies in the class $\alpha$.
\end{remark}

\begin{remark}[G\"{u}rel's theorem]\label{rem:Gurel}
    Proposition \ref{prop:existence_of_non-contractible_orbits} reflects the same filtered-Floer mechanism that underlies existence results for infinitely many non-contractible periodic orbits, for which see \cite{Gu13,GG16,Or17,Or19,Su21}. In particular, \cite[Theorem 1.1]{Gu13} proves the following. Let $H$ be a Hamiltonian on a closed manifold with atoroidal symplectic form, and suppose that $H$ has a non-degenerate one-periodic orbit $x$ whose homology class $[x]$ is non-zero in $H_1(M;\ZZ)/\mathrm{Tor}$ and that the set of one-periodic orbits representing $[x]$ is finite. Then, for every sufficiently large prime $p$, $H$ has a simple periodic orbit in the homology class $p[x]$ whose period is either $p$ or the next prime $p'$.

    A derivation of that prime-period theorem requires additional iteration input: admissible iterations preserve the relevant local Floer homology, the action of the iterated orbit is homogeneous, and one needs action-gap and prime-gap estimates (for the latter, see \cite{BHP01}). These facts do not follow merely from the rectangle barcode or from Proposition \ref{prop:existence_of_non-contractible_orbits}. We therefore do not identify the side lengths of an individual rectangle with homogeneous quantities under iteration without a separate argument. See \cite{Gu13} for the complete iteration proof.
\end{remark}

\section*{Acknowledgments}

The authors are deeply grateful to the anonymous referee for exceptionally valuable contributions to this paper, for permission to include the results and proofs described in the Attribution subsection, and for careful reports that substantially improved the exposition.
We also thank Qi Feng and Vadim Lebovici, whose careful readings of earlier versions uncovered an error concerning rectangles of type (R), and who independently showed that no such rectangle occurs. Feng gave a singular-value-decomposition approach, while Lebovici pointed out the relation to homology pyramids in Remark \ref{rem:pyramid} and the relevance of the decomposition theorems of \cite{CGL21} and \cite{UZ16}, and both observed the connection between the spectral spread, the finite bar lengths, and the boundary depth. We are grateful to them for generously sharing these insights.
We also thank Tomohiro Asano for fruitful discussions.

\section*{Declarations}

\textit{Funding.} R.~Orita was supported by JSPS KAKENHI Grant Number 21K13787.

\textit{Competing interests.} The authors have no competing interests to declare that are relevant to the content of this article.

\textit{Data availability.} Data sharing is not applicable to this article, as no datasets were generated or analysed.

\bibliographystyle{alpha}
\bibliography{orita_bibtex}

\end{document}